\documentclass[final,11pt,onecolumn,a4]{IEEEtran}

\usepackage{bm,graphicx,epsfig,amsmath,amssymb,mathrsfs}
\usepackage{cite}

\newlength{\intwidth}
\DeclareRobustCommand{\fpint}[2]
  {\mathop{%
     \text{%
       \settowidth{\intwidth}{$\int$}%
       \makebox[0pt][l]{\makebox[\intwidth]{$-$}}%
       $\int_{#1}^{#2}$}}}

\newtheorem{theorem}{Theorem}

\newtheorem{corollary}{Corollary}

\newtheorem{definition}{Definition}[section]

\newtheorem{lemma}{Lemma}

\newtheorem{proposition}{Proposition}

\newcommand{\E}{\ensuremath{\mathrm{E}}}
\newcommand{\var}{\ensuremath{\mathrm{var}}}
\newcommand{\cov}{\ensuremath{\mathrm{cov}}}

\newcommand{\f}{{\bf f}}
\newcommand{\x}{{\bf x}}
\newcommand{\kk}{{\bf k}}
\newcommand{\y}{{\bf y}}

\begin{document}

\title{ Hyperanalytic Denoising}

\author{Sofia~C.~Olhede \\
Department of Mathematics, Imperial College London, SW7 2AZ UK% <-this % stops a space
\thanks{
This work was supported by an EPSRC grant.}%<-this % stops a space
\thanks{S. Olhede is with the Department of Mathematics,
Imperial College London, SW7 2AZ, London, UK (s.olhede@imperial.ac.uk). Tel:
+44 (0) 20 7594 8568, Fax: +44 (0) 20 7594 8517.}}

\markboth{Statistics Section Technical Report TR-06-01}{Olhede:
Hyperanalytic Denoising}
\maketitle

\begin{abstract}
A new thresholding strategy for the estimation of a
deterministic
image
immersed in noise is introduced. The
threshold is combined with a wavelet decomposition, where the
wavelet coefficient of the image at any fixed value of the decomposition
index is estimated, via thresholding the observed coefficient depending on
the value of both the magnitude of the observed coefficient as well as the magnitudes of coefficients
of a set of additional images calculated from the observed image. The additional
set of images is chosen so that
the wavelet transforms of the full set of images have suitable deterministic
and
joint stochastic properties at a fixed scale and position index.
Two different sets of additional images are suggested. The behaviour of the
threshold criterion
for a purely noisy image is investigated and 
a universal threshold is determined. The properties of the threshold
for some typical deterministic signal structures are also given.
The risk of an individual coefficient is determined, and calculated explicitly
when the universal threshold is used, and some typical deterministic signal
structures. The
method is implemented on several examples and the theoretical risk reductions
substantiated.
\end{abstract}

%EDICTS: MRP-WAVL (Wavelets), RST-DNOI (Denoising).                  
\section{Introduction \label{Introduction}}
\PARstart{T}{his} paper treats the problem of estimating an unknown deterministic
image
immersed in noise. The proposed estimation procedure will be based on a separable
wavelet decomposition of the observed image that is augmented by a set of
wavelet decompositions of additional images calculated from the
observed image. The wavelet coefficients
of the full set of images at any fixed value of the scale and position
are 
used to estimate the wavelet transform coefficient of the deterministic image at the given scale and position. The transform is then inverted
and the spatial domain image
estimated. In 1-D signal estimation Donoho and Johnstone \cite{Donoho1994,Donoho1995} first
proposed
estimation of a noisy deterministic signal using the wavelet transform. 
The success of such decomposition based methods mainly relies on the deterministic and stochastic properties
of the observed or noisy decomposition coefficients at any fixed index value.
In the simplest form the estimation procedure
roughly corresponds to separating `clean' and `noisy' coefficients into two subsets, where the `noisy' coefficients
are eliminated or subjected to some form of shrinkage \cite{Stein}.
Often
each coefficient is estimated separately at any given index value and for example the procedure
may correspond to eliminating the coefficients whose magnitudes do not exceed
a given threshold.
A possible choice of threshold is the universal threshold, constant across
coefficients, that for large sample sizes gives a risk close to that given
by using an `oracle,' {\em i.e.} knowing whether a coefficient should be
eliminated
or retained \cite{Donoho1994}. A slightly different definition is
given to the universal threshold by the authors of \cite{downie}, that we
shall use. If the
decomposition is highly
compressed hard thresholding combined with the universal threshold will achieve very
good estimation in terms of low mean square error.

Naturally, to achieve optimal compression for locally simpler structures,
such as 1-D behaviour embedded in a 2-D image, whilst still
being able to represent varied signal structure in 2-D,
the decomposition algorithm becomes more complicated. If a very simple decomposition method is used, then determining
the statistical properties of the observed coefficients is easily done, and
the decomposition
can be found without major computational expense. The draw-back is that in
general the
mean square error of the estimation will, unless the estimation procedure
is more complicated, with a very simply decomposition often
increase due to lack of compression. Hence a trade-off
must be found between the choice of decomposition and the appropriate treatment of the coefficients
of the observed
signal to form estimates.
This compromise
will naturally vary with the assumptions placed on the observed signal. As the variational
structure in 2-D is a great deal richer than in 1-D, many different methods
have been developed to achieve optimal compression for given image structures, and for
particularly successful
examples see work by Starck {\em et al.} on curvelets \cite{starck}, work by Donoho
on wedgelets \cite{donoho} or work by le Pennec and Mallat on bandelets \cite{pennec}.
An important feature of all these decompositions is the representation
of
an image in terms of coefficients associated with a given spatial position,
and length scale. The coefficients are considered `local' to such positions
and scales.

\begin{figure*}[t]
\centerline{
\includegraphics[scale=0.97]{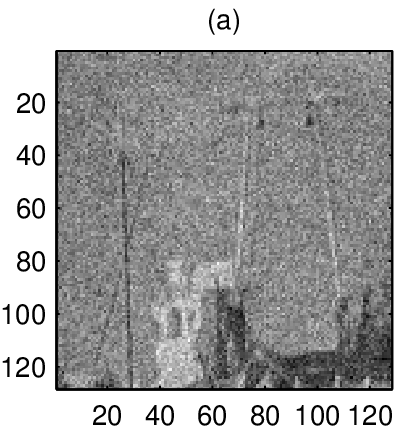}
\includegraphics[scale=0.97]{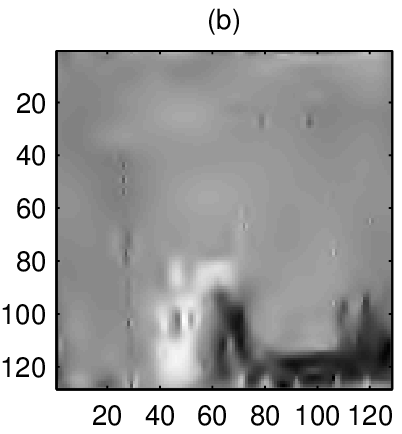}
\includegraphics[scale=0.97]{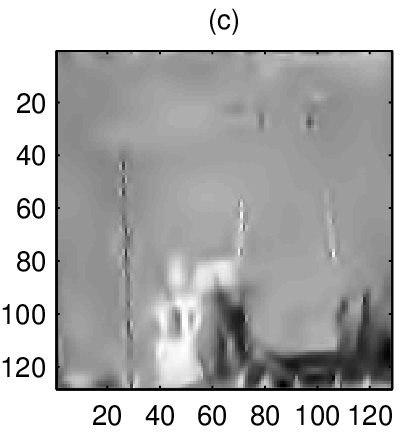}
\includegraphics[scale=0.97]{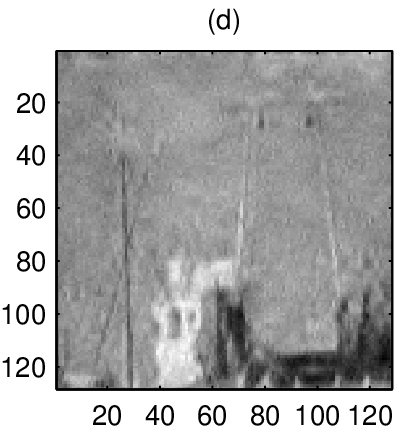}}
\centerline{
\includegraphics[scale=0.95]{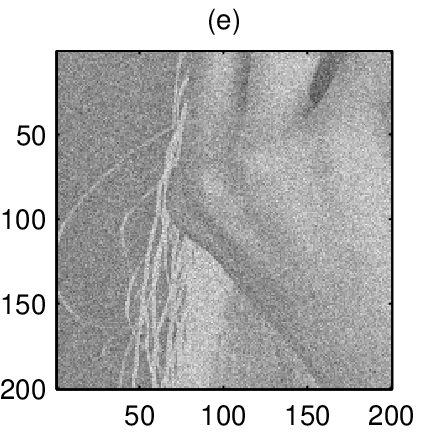}
\includegraphics[scale=0.95]{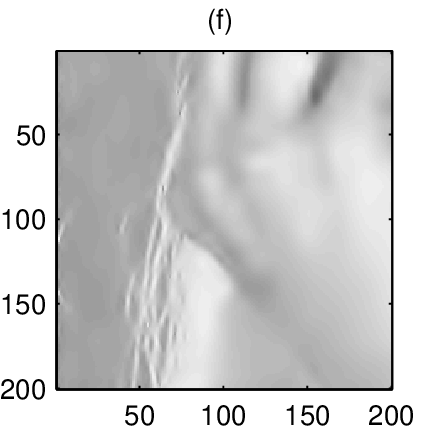}
\includegraphics[scale=0.95]{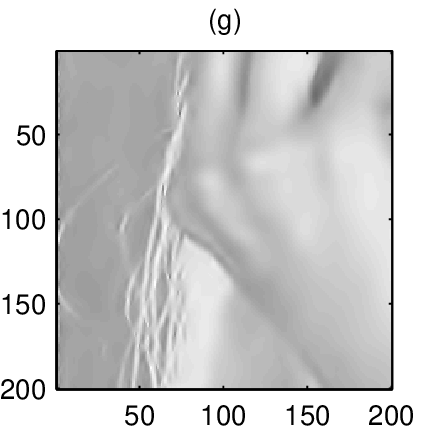}
\includegraphics[scale=0.95]{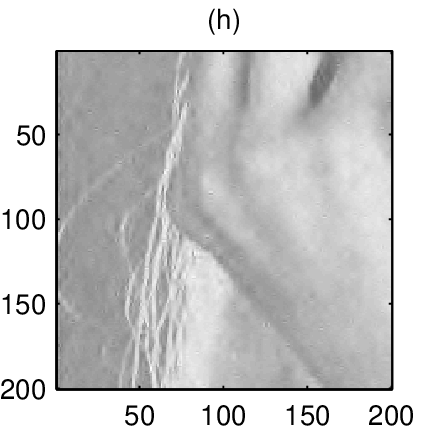}
}
\caption{\label{boat}
Row 1: estimating a section of the noisy boat image (SNR=5.56), with the noisy image (a), the usual hard thresholding
estimate (b), the hypercomplex estimate (c) and the HMM estimate (d). The
hypercomplex estimate in (c) achieves greater continuity, see for example the mast, than that achieved by (b). The HMM estimate has more noise left in the reconstruction, a feature also observed
by Starck {\em et al.} \cite{starck}[p.~680].
Row 2: a section of one of the bands of the noisy Tiffany image (SNR=8), with the noisy image (a), the usual hard thresholding
estimate (b), the hypercomplex estimate (c) and the HMM estimate (d). Observe
the curved loose strand of hair, and that the HMM method reconstructs
the background with some noisy artifacts.
}
\end{figure*}

Methods that achieve a substantial degree of compression
can afford to
treat each decomposition
coefficient individually and without a great deal of sophistication. To achieve better estimation of coefficients for a simplistic decomposition
method treating decomposition coefficients
{\em simultaneously} may also give improved estimation. This may correspond to full likelihood
based or similar methods such as those proposed by Jansen and Bultheel \cite{jansen2001}, note also work by Johnstone and Silverman where the local
sparsity of decomposition coefficients is discussed 
\cite{johnstonesilverman2004}, as well as usage of specific known coefficient structure in the decomposition of a deterministic
signal across
coefficients: see for example Cai and Silverman
\cite{Cai}, Dragotti
and Vetterli \cite{dragotti},
Pi{\v z}urica {\em et al.} \cite{pizurica}, Crouse {\em et al.} \cite{crouse},
Fryzlewicz \cite{fryz} and Olhede and Walden \cite{olhede2004}. By modelling
continuity across coefficients in terms of their local index, estimates with
a
reduced mean square error may be obtained, that frequently correspond to
better visual reconstructions of the image.
If the method of estimating the decomposition coefficients
is not very complicated but still captures continuity across the decomposition
index well, then we may choose
a decomposition of the image that is not optimal in terms of compression
of the deterministic image, but that is computationally cheap to implement, and where the estimation of the
decomposition coefficients may be treated carefully. We may then
achieve a reduced mean square error in the estimation of
the image at a low computational expense. The purpose of this paper
will {\em not} be to develop an optimal decomposition algorithm,
but instead to improve the estimation of the decomposition coefficients,
without complicating the procedure substantially. We shall base our image estimate on the 2-D separable wavelet transform
coefficients, extending 1-D methods of utilizing coefficient structure
to 2-D. We shall make the developments precise by discussing the risk of
a given estimated decomposition coefficient.

In 1-D Dragotti and Vetterli \cite{dragotti} explicitly model signal
structure as
a polynomial function plus some discontinuities and jointly estimate the full set
of coefficients corresponding to a discontinuity. Pi{\v z}urica {\em et al.}
pool information regarding joint structure in 2-D across coefficients via estimating the local Lipschitz constant and this information is used to estimate
the probability that a wavelet coefficient contains contributions from a
signal. Crouse {\em et al.} \cite{crouse}[p.~887] model signal
presence across coefficients in 1-D in terms of clustering and persistence, i.e.
if a coefficient is non-null at a given decomposition index, then coefficients
that are ``close'' to this index are also non-null. A similar strategy
is adopted in 1-D by Cai and Silverman \cite{Cai}, whilst Fryzlewicz
\cite{fryz} considers the magnitude of any additional arbitrary coefficient
when estimating a given coefficient.
Fryzlewicz established the risk of this strategy, determined from the mean and
covariance matrix of the two coefficients. Fryzlewicz's treatment is very general
and instructive. 

In a similar spirit to some of the aforementioned methods in 1-D Olhede and Walden \cite{olhede2004} considered the thresholding
of an individual wavelet coefficient based on the magnitude of the observed
coefficient and the magnitude of the decomposition coefficient of the Hilbert
Transform (HT) of the signal, denoting this method `analytic' denoising. As both the HT and the wavelet transform are linear the strategy can be viewed either as constructing a second out-of-phase replication of
the original signal and finding its local decomposition, or as forming a
weighted average of coefficients of the same scale that are nearby in time
and using this magnitude to determine if there is local signal presence.
The latter strategy is similar in spirit to block thresholding,
but instead of using a local magnitude calculated from an average
of squared adjacent wavelet coefficients at the
same scale in the thresholding procedure, the square of a weighted average of coefficients with a weighting of $O\left(1/(2^j k)\right)$ is used. For `analytic' thresholding to work well the wavelet
coefficient of the HT of the deterministic signal must be large when the wavelet transform
of the signal should not be estimated by zero, and the distribution of the
wavelet coefficient of the HT of the noise needs to be jointly determined
with the wavelet coefficient of the noise at the same scale and time. As the HT can be considered to have
the same time-frequency structure as the original signal the wavelet transform
of the signal and the HT of the signal should be large concurrently. 
Olhede and Walden \cite{olhede2004} determined that the wavelet transform of the
noise and its HT were approximately uncorrelated at a fixed time and scale, and supplied an appropriate
universal threshold for `analytic' denoising. Figure \ref{exx2}(a) shows the risk of a given coefficient using
`analytic' denoising, based on the wavelet transform of  the signal and its HT taking the same magnitude.
The figure verifies that in the case of equivalent means the risk of an `analytic' hard thresholded coefficient estimate with a universal threshold is less than that of a hard thresholded estimate with a universal
threshold. Thus
improvements to standard denoising in 1-D can be obtained by implementing
this
procedure.

We seek to extend `analytic' thresholding to
2-D, and this will in general require defining additional images, serving the same purpose as the HT did in 1-D. The HT was useful
in 1-D, as it has the same
time-frequency structure as the original signal, something we discuss in section IIB, and also the joint statistical
properties of the wavelet transform of noise and its HT at a fixed scale and position was easily determined.
In 2-D there are many possible
extensions to the HT, where in each case more than a single
additional component is defined. We refer to such components as {\em quadrature
components,} that are introduced and discussed in sections IIB and IIC.
There is more than one extension because variation in the image can
either be locally uni-directional, and associated with a given direction, or occurring
in several directions simultaneously (see Olhede and Metikas \cite{Metikas2}
for a more complete discussion of this topic). We investigate the usage of two possible HTs: the Riesz Transforms (RTs, section IID) of the image
or the tensor
products of the HT in 1-D with the identity filter, denoted the HyperComplex
transforms (HCTs, section IIE). We define the local magnitude of the wavelet
coefficients from the quadrature components (section IIIA), and propose a
threshold criterion to estimate the wavelet coefficients of the image. Once
the wavelet transform is inverted this yields an estimate of the image, and
this method is denoted hyperanalytic denoising.

We discuss the properties of the local magnitude
for stylized image structure: {\em i.e.} the behaviour of the threshold
criterion for oscillatory structures and edges (section IIIB). We discuss
the choice of threshold, and an appropriate universal threshold for correlated
threshold criteria (section IIIC).
We determine
the approximate distribution of the decomposition of the Riesz and Hypercomplex
components
of noise alone
at a fixed value of the indexing (section IIID), and this allows us to determine
universal
thresholds for both the RT and HCT based methods (section IIIE).
We calculate the approximate risk associated with the two different
thresholding strategies with a given threshold (section IIIF), and discuss
the value of the risk of the different procedures for certain scenarios.

We implement the procedure on several examples (section IV), and compare
results with the Hidden
Markov Model method (HMM). We observe that a reduced mean square error is
obtained from using the proposed image
denoising strategies, and discernable improvements in the visual reconstructions.
Hyperanalytic denoising is thus shown to give a simple and
computationally
competitive method of improving existing denoising strategies.

\section{Image Model}
\subsection{Image Structure}
We model the observed image $\left[Y_{\x}\right]_{\x}$ for 
$\x=\left[x_1,x_2\right]^T\in D,$ where $D=\left[0,N-1\right]^2,$ and $\Delta x$ denotes the sampling period
via:
\begin{equation}
Y_{\x}=q(x_1 \Delta x,x_2 \Delta x)+\epsilon_{\x},\quad \x\in D
\label{model}.
\end{equation}
We collect the observed image in a matrix $\bm{Y}=
\left[Y_{\x}\right]_{\x\in D},$ and similarly define $\bm{q}=\left[q_{\x}
\right]_{\x\in D}=
\left[q(x_1 \Delta x,x_2 \Delta x)\right]^T_{\x\in D},$ as well as $\bm{\epsilon}=
\left[\epsilon_{\x}\right]_{\x\in D} .$
The noise is modelled by
$\bm{\epsilon}_{\x}\sim N\left(0,\sigma^2\right),$ where $\sim$
denotes distributed as, and 
$
{\mathrm{Cov}}\left(\bm{\epsilon}_{\x},\bm{\epsilon}_{\y}\right)=
\sigma^2\delta_{\x,\y},\quad \x,\;\y\in D,$
i.e. the noise is Gaussian, uncorrelated and isotropic. A decomposition
of the image in terms of a wavelet basis \cite{mallat} is formed via
\begin{equation}
%\nonumber
q(\x \Delta x)=
\sum_{j,\kk} W_{j,1,\kk}^{(q)} \psi_{j,1,\kk}(\x)
+\sum_{j,\kk} W_{j,2,\kk}^{(q)} \psi_{j,2,\kk}(\x)+
\sum_{j,\kk} W_{j,3,\kk}^{(q)} \psi_{j,3,\kk}(\x)+
\sum_{\kk} V_{k_1,k_2}^{(q)} \phi_{J,\kk}(\x),
\label{inverting}
\end{equation}
where  $
\psi_{j,1,\kk}(\x),$
$\psi_{j,2,\kk}(\x),$ 
$\psi_{j,3,\kk}(\x),$ and $
\phi_{j,\kk}(\x)$ are the tensor products of functions
$\psi_{j,k}(x)$ and $\phi_{j,k}(x),$ respectively.
$V_{J,\kk}^{(q)}\equiv W_{J,4,\kk}^{(q)}$ is then associated with smooth behaviour in the image
$q(\x)$ in the variables $x_1$ and
$x_2$,
$W_{j,2,\kk}^{(q)}$ is associated with smooth behaviour in $x_2$ and
rapid variation in $x_1,$ {\em etc,} where $u=1,\dots, 4$ denotes the tensor
product index. $j$ is associated with scale $2^j,$ where $1\le j\le J_0={\mathrm{lg}}(N),$
whilst $\kk$ is associated with a spatial localisation
in the plane. If an image with $N^2$ coefficients is observed, then for any
fixed value $j,$ $0\le k_l \le N_j-1,\;l=1,\;2,$ where $N_j=N/2^j.$ 
For simplicity collect the indices in a vector-valued index
of $\bm{\xi}=\left[j,u,\kk\right]^T.$ The full set of coefficients 
$\left[W_{j,u,\kk}^{(q)}\right]$ is the Discrete Wavelet Transform (DWT)
of $q(\cdot).$

The DWT is usually implemented by
repeated filtering of the
observed
signal with two special filters, the scaling filter $\left\{g_l:\;l=0,\dots,L-1\right\}$
and the wavelet filter $\left\{h_l:\;l=0,\dots,L-1\right\},$ in both spatial
directions separately. We initialise the transform by equating the image with
the finest scale representation of
the image, i.e. $V_{0,x_1,x_2}^{(q)}\equiv q_{x_1,x_2}.$ 
The transform at index $\bm{\xi}$ can also be implemented using a single filter
$h_{j,l}.$ 
The decomposition is halted at level $j=J\le J_0={\mathrm{lg}}(N),$ and
the scaling coefficients $\left\{V_{J,\kk}^{(q)}\right\}$ are determined
at this level to complete the representation. Hence for $j<J$ only
$\left[W_{j,u,\kk}\right]$ for $u=1,2,3$ are calculated.
For more details on the DWT, see for example Percival
\& Walden \cite{percival}, whilst a good exposition of image decompositions
can be
found in Mallat \cite{mallat}. Having observed $Y_{\x}$ rather than
$q_{\x}$ we calculate the DWT coefficients $W_{\bm{\xi}}^{(Y)},$
and threshold these to obtain an estimate of $W_{\bm{\xi}}^{(q)},$ denoted
$\widehat{W}_{\bm{\xi}}^{(q)}.$
Wavelets will compress images of sufficient regularity, a statement that
can be made precise in terms of Besov spaces,
but for some locally simple image structures, a more compressed representation can be made
\cite{starck,pennec}. Hence for images containing say edges the deterministic
image energy in the DWT will be spread over more coefficients than strictly necessary,
and as the magnitude of the affected coefficients will be less than
the coefficients representing the same structure in a more compressed
alternative decomposition it is important that the estimation procedure does
not fail to retain signal generated coefficients.

\subsection{Quadrature Components}
In one version of the 1-D estimation
algorithms suggested by Cai and Silverman \cite{Cai}, the coefficient at
scale $j$ and position $k$ was estimated using a shrinkage rule depending on the combined
magnitude of the observed coefficient at $[j,k]$ and the magnitude of the immediate time-neighbours
at the same scale, i.e. at $[j,l]$ for $l\neq k.$ This procedure will perform well if a signal contribution present at the
$[j,k]$ index exhibits clustering in adjacent coefficients,
%, see \cite{crouse}[p.~887], 
{\em i.e.} the wavelet coefficients
will have large means at $[j,l],$ and the noise is uncorrelated
over $l\neq k.$ The scale adjacent coefficients
at a given time point have been used to improve estimation \cite{fryz,crouse}, and Olhede and Walden
\cite{olhede2004} used the wavelet decomposition of the HT of the observed image to this purpose. We seek to generalise the method in \cite{olhede2004}
to 2-D, and discuss some of
its properties, before proceeding to do so.

To simplify the discussion of the HT, let the Fourier Transform (FT) of a
$d$ dimension signal
$q(\x)$ be denoted by:
\[Q(\f)=\int_{{\mathbb{R}}^d} q(\x) e^{2\bm{i}\pi \f^T\x}\;d^d\x=
\left|Q(\f)\right|e^{-2\bm{i}\pi \varphi_q(\f)},\]
this defining the magnitude ($\left|Q(\f)\right|$) and phase 
($\varphi_q(\f)$)
of $q(\x)$ in the Fourier domain.
Given a 1-D signal $q(x)$ the HT in both the time and frequency
domain are defined by:
\begin{equation}
\label{hilberttrans}
{\cal H}q(x)=\frac{1}{\pi}\fpint{-\infty}{\infty}
\frac{q(y)}{x-y}\;dy,\quad
\left({\cal H}Q\right)(f)=(-\bm{i})Q\left(f\right){\mathrm{sgn}}(f),
\end{equation}
and the transform can be approximated suitably for discrete implementation (see \cite{olhede2004}). 
Usually $q(x)$ and ${\cal H}q(x)$ are collected into a complex-valued 
representation, denoted
the analytic signal, given by $q^{(+)}(x)=q(x)+\bm{i}{\cal H}q(x).$
If $q(x)=\cos(2\pi f^{\prime} x)$ then $q^{(+)}(x)=\exp(2\pi\bm{i} f^{\prime} x),$ but sometimes too much emphasis is put on this description of oscillatory
signals, to the extent where the HT is almost discounted in usage when the
observed signal is not oscillatory.
Even if $q(x)$
does {\em not} correspond to an oscillation, the HT can be considered to
enjoy certain properties, such as: i)
${\cal H}q(x)$ is orthogonal to $q(x),$ {\em i.e.} $\int {\cal H}q^*(x)q(x)\;dx=0,$
ii) the magnitude of the HT of $q(x)$
at any given frequency $f\neq 0$ is identical to that of
the original signal, {\em i.e.} $\left| Q(f)\right|^2=\left|{\cal H}Q(f)\right|^2,$
iii) the HT is linear in the signal,
and iv) the HT of a signal can be considered as having the
same time-frequency signature as the original signal. i-iii) immediately
follow from equation (\ref{hilberttrans}), and ensure
that the distribution of the DWT of the HT of white noise at a given value
of $[j,k]$ is asymptotically identical to that of the DWT of
the original noise, and the two wavelet coefficients are approximately uncorrelated
\cite{olhede2004}.
The fourth property merits
some further discussion. Clearly the notion that a signal and its HT have the same time-frequency structure is accepted in signal processing,
as the analytic signal, rather than the real signal, is used to construct
time-frequency representations of a real signal.
For example usage of the Wigner-Ville distribution rather than the Wigner distribution
is generally advocated \cite{boashash}. As may be noted from equation (\ref{hilberttrans})
${\cal{H}}(q)(x)$
has at all frequencies $f\neq 0$ exactly the same frequency support as $q(x),$
whilst the spatial support of $q(x),$
has been spread out by the convolution with $1/(\pi x).$ 

The HT
is usually interpreted as a phase-shift of $\pi/2$ to signal
$q(x).$ Note that we may write:
\begin{equation}
\label{phase}q(x)=2\int_{0}^{\infty} 
\left|Q(f)\right|\cos(2\pi( fx-\varphi_q(f)))\;df,\;
{\cal H}q(x)=2\int_{0}^{\infty} 
\left|Q(f)\right|\sin(2\pi( fx-\varphi_q(f)))\;df.
\end{equation}
Hence the same magnitude of $\left|Q(f)\right|$ is assigned to each frequency
$f,$ and the contribution previously associated with $\cos(2\pi (\varphi_q(f)- fx))$ is now shifted in cycle or phase by $\pi/2.$ Thus in some sense, we are recovering
the same signal, as the frequency description is the same, but there has
been a very slight shift in time alignment. 
Thus ${\cal H}q(x)$ should have roughly
the same
time-frequency support as $q(x).$ This implies that the DWT
coefficient of $q(x)$ should have about the same magnitude as ${\cal H}q(x),$
as the DWT forms a time-frequency decomposition of a given
signal.
The DWT is compact in time, and we wish to encourage using time information
in nearby locations when estimating a given coefficient in analogue with
Cai \& Silverman. Of course once the coefficient has been estimated,
the estimate of the signal will be based on the thresholded wavelet coefficients
of the observed signal alone, and thus discontinuities can still be reconstructed.

Given the nice deterministic and stochastic properties of `analytic' denoising,
it is not strange that we seek to generalise the concept to 2-D. A first step in this procedure is the definition of
linear transformations of the image that will serve the same purpose as the
HT did in 1-D. The HT and the signal formed a natural representation
in terms of the `analytic' signal, where the real and imaginary components
were phase shifted versions of each other, or we may denote the latter two signals
as being `in quadrature,' {\em i.e.} as representing out-of-phase replications
of the same structure. Their magnitude squares represent the local
presence of the signal well but we stress that even if the signal is not oscillatory, the
interpretation of the HT as having roughly the same time-frequency support
still rests on the above arguments.
We shall denote the signal and its HT
as quadrature components, and will define the 2-D generalization of this
two-component signal collection.

\begin{definition}[Quadrature Components]
We denote by {\em Quadrature Components} of $q(\x)$ any set of images
$\left\{\breve{q}^{(s,l)}(\x)\right\}_{l=1}^L,$
where $s$ denotes the specific transform used in the construction of the
components that satisfy:
\begin{enumerate}
\item \label{cond1} each $\breve{q}^{(s,l)}(\x)$ is orthogonal (`out of phase') to the
original signal
$q(\x)\equiv\breve{q}^{(s,0)}(\x),$ or \\
$\int_{-\infty}^{\infty} \int_{-\infty}^{\infty}q^{*}(\x) \breve{q}^{(s,l)}(\x)\;d^2\x=0,$
and for all separable 
$q_S(\x)=q_1(x_1)x_2(x_2),$ also\\
$\int_{-\infty}^{\infty} q^{*}_S(\x) \breve{q}_S^{(s,l)}(\x)\;dx_l=0,$ 
\item \label{cond2} the combined energy assigned to each frequency $\f$ from the full set
of quadrature components at all points of $\f\in\mathbb{R}^2$ except for
a finite set of frequencies satisfies the equation
$
\sum_{l=1}^L \left|\breve{Q}^{(s,l)}(\f)\right|^2=C_L^{(s)}\left|Q(\f)\right|^2,
$
where $0<C_L^{(s)}<\infty$ is constant and,
\item \label{cond3} each $\breve{q}^{(s,l)}(\x),$ for $l=1,\dots,L,$ is constructed by a
linear transformation of $q(\x),$
\item \label{cond4} the space and spatial frequency support of $\left[\breve{q}^{(s,1)}(\x),
\dots,\breve{q}^{(s,L)}\right],$ for $l=1,\dots,L$ is similar to that of
$q(\x).$
\end{enumerate}
\end{definition}

We form the DWT of all $L+1$ images, $\left\{\breve{q}_{\x}^{(s,l)}
\right\}_{l=0}^L,$ and define:
\begin{equation}
\label{wavcoeff7}
\breve{W}_{\bm{\xi}}^{(q,s,l)}=W_{\bm{\xi}}^{(\breve{q}^{(s,l)})},\;l=0,\dots,L.
\end{equation}
The linear operator that constructs object $\breve{q}^{(s,l)}(\x)$ from $q(\x)$
will be denoted ${\cal V}^{(s,l)}$ and the transformation is implemented
in the spatial domain using the kernel $v^{(s,l)}(\x)$ that once the integral
is approximated using a Riemann sum is replaced by a
linear filter $v_{D,\x}^{(s,l)}.$ The FT of 
$v^{(s,l)}(\x)$ is denoted $V^{(s,l)}(\f)$ whilst taking the FT
of $v_{D,\x}^{(s,l)}$ yields the object $V_D^{(s,l)}(\f).$ The discrete
implementation of the calculation of the quadrature components is outlined
in Appendix A.

\subsection{Stochastic Properties of Decomposed Quadrature Components}
We establish the stochastic properties of the wavelet decomposition of
noise alone, and for this purpose define at a fixed value of $\bm{\xi}:$
\begin{equation}
\label{noisecollection}
{\mathbf{n}}^{(s,u)}=
\left[\breve{W}_{\bm{\xi}}^{(\epsilon,s,0)},\dots,
\breve{W}_{\bm{\xi}}^{(\epsilon,s,L)}\right]^T,\;s=r,\;h.
\end{equation}

\begin{proposition}[Energy of Quadrature Components \label{proppy1}]
At a fixed index value $\bm{\xi}=[j,u,\kk]^T$ the total energy of the DWT
of the quadrature
components of white noise with variance $\sigma^2$ is given by:
\begin{equation}
\E\left(\sum_{l=1}^{L} 
n^{(s,u)2}_l\right)=C_{L}^{(s)} \sigma^2+O\left(1/N\right),\;u=1,2,3,4.
\end{equation}
\end{proposition}
\begin{proof}
See appendix C. The error term follows from the Riemann approximation
to the integral.
\end{proof}

\begin{proposition}[Covariance of Transforms of the Signal and Its Quadrature Components \label{proppy2}]
At a fixed index value $\bm{\xi}=[j,u,\kk]^T$ the covariance of the DWT of white
noise, and the DWT of any of the quadrature components of the white
noise
if of order $O(1/N).$
\end{proposition}
\begin{proof}
See appendix C. The error term follows from the Riemann approximation
to the integral.
\end{proof}

Thus at any given value of $\bm{\xi}$ the DWTs of $\epsilon_{\x}$
and $\left\{\breve{\epsilon}_{\x}^{(s,l)}\right\}$ are approximately uncorrelated,
and the combined energy of the DWTs of 
$\left\{\breve{\epsilon}_{\x}^{(s,l)}\right\}$ is a multiplicative constant
of the energy of the DWT of  $\epsilon_{\x}.$ Thus the squared magnitudes
of these objects have a tractable joint distribution. Condition \ref{cond4} ensures that we may assume that
the mean of the DWT of the observed image will be simultaneously large to the mean of
the DWTs of the quadrature components of the observed image at a given value of the index $\bm{\xi}.$  Of course whilst the general definition of `quadrature components' may then seem
justifiable, this does not guarantee the existence of such objects. We shall
give two different specific examples of quadrature components based on extending
the
HT to 2-D, and discuss
their properties. 
We base the set of quadrature components on hyperanalytic functions,
see \cite{Metikas2}.

\subsection{The Riesz Transforms \label{riessy}}
The Riesz Transforms (RTs) have been used in combination with
the wavelet transform by Metikas and Olhede \cite{Metikas,Metikas2}.
Denote the convolution of two functions $q_1(\x)$ and $q_2(\x)$ by
$(q_1 \ast \ast q_2)(\x)=\int \int_{{\mathbb{R}}^2} q_1(\y) q_2(\x-\y)\;d^2\y.$
The RTs of $q(\cdot),$ denoted $\breve{q}^{(r,1)}(\x)$ and $\breve{q}^{(r,2)}(\x),$ are obtained by convolving $q(\cdot)$ with the Riesz kernels $v^{(r,l)}(\x),$
given in terms of $x=\sqrt{x_1^2+x_2^2}$ and $f=\sqrt{f_1^2+f_2^2}$ by:
\begin{equation}
v^{(r,l)}(\x)=\frac{x_l}{2\pi x^3} ,\quad V^{(r,l)}(\f)=-\bm{i}
\frac{f_l}{f},\;l=1,\;2,\quad
\breve{q}^{(r,l)}(\x)=\left(v^{(r,l)}\ast \ast q \right)(\x),\;l=1,\;2.
\end{equation}
The RTs satisfy the conditions of quadrature components, see for example
\cite{Metikas2}[p.~15--16]. Given the RTs combine to have the same norm as
$q(\x),$ $C_2^{(r)}=1.$ As in Olhede \& Metikas 
\cite{Metikas2} we argue that
if unidirectional structure only is present in
the image with orientation $\nu,$ i.e. the image admits the representation for $\mathbf{n}=\left[\cos(\nu)\quad \sin(\nu)\right]^T$ and $\nu\in(0,\pi)$
of:
\begin{equation}
\label{unidir}
q_U(\x)=\int_{0}^{\infty} G_U(f) \cos(2\pi f\mathbf{n}^T\x)\;df,
\end{equation}
then the interpretation of the RTs is simplified. We use polar coordinates
and set $\f=f\left[\cos(\phi)\quad \sin(\phi)\right].$
Then the Fourier transform of $q_U(\x)$ is $Q_U(\f)=\frac{G_U(f)}{2}\left[
\delta(\phi-\nu)+\delta(\phi+\nu) \right]$
and we find:
\begin{equation}
\label{riesz1}
\left[\breve{q}_U^{(r,1)}(\x)\;\breve{q}_U^{(r,2)}(\x)\right]=
\left[\cos(\nu)\;\sin(\nu)\right]  \int_{0}^{\infty} G_U(f)
\sin(2\pi f\mathbf{n}^T\x)\;df.
\end{equation}
Thus the two quadrature components represent the same 1-D directional variation as $q_U(\x),$ with
the same directionality as $q_U(\x,)$ but where the variations in direction
$\nu$ have been shifted in phase and multiplied by a constant factor. Thus (informally) for unidirectional variation the Riesz transforms have the same spatial and
spatial frequency support as the original signal.
Note that we are {\em not} assuming that $q(\x)$ is periodic or oscillatory.

\subsection{The Hypercomplex Transforms \label{hyperc}}
A second set of 2-D HTs are
the HyperComplex Transforms (HCTs), defined as tensor products of the identity
transform and the
HTs. By
Olhede and Metikas \cite{Metikas2}[p.~12--13], it is shown that the hypercomplex
transforms give a valid set of quadrature components, and note $C_3^{(h)}=3.$ Denote the partial HT \cite{hahn1996}
in
direction $x_l$ by ${\cal H}_l.$ Three additional quadrature components are defined by:
\begin{eqnarray}
\breve{q}^{(h,1)}(\x)={\cal H}_1\left\{q\right\}(\x),\quad
\breve{q}^{(h,2)}(\x)={\cal H}_2\left\{q\right\}(\x),\quad
\breve{q}^{(h,3)}(\x)={\cal H}_2{\cal H}_1\left\{q\right\}(\x).
\end{eqnarray} 
If the image is naturally expressed as separable in the frame of reference
the three HCTs of $q_S(\cdot)$ are by trivial extension of equation (\ref{phase}), the same signal shifted in phase
in the two axes. Of course the
purpose of this paper will be to alleviate problems (see for example
Starck {\em et al.} \cite{starck}[p.~671]
) when estimating nonseparable
images
based on coefficients calculated in a separable decomposition 
whose
energy spread across more coefficients than strictly necessary. Assume $q(\x)$ is non-separable then 
define its Partial FT (PFT) in direction $x_1$ by:
$Q_1(f_1,x_2)=\left|Q_1(f_1,x_2)\right|e^{-
2\pi\varphi_1(f_1,x_2)}=
\int_{-\infty}^{\infty} q(\x) e^{-2\bm{i}\pi f_1 x_1} \;dx_1,$
so
\begin{eqnarray}
q(\x)&=&2\int_{0}^{\infty} \left|Q_1(f_1,x_2)\right| \cos(2\pi (
f_1 x_1-\varphi_1(f_1,x_2)) \;df_1\\
{\cal H}_1 q(\x)&=&2\int_{0}^{\infty} \left|Q_1(f_1,x_2)\right| \cos(2\pi
(
f_1 x_1-\varphi_1(f_1,x_2)-\pi/2) \;df_1.
\end{eqnarray}
Thus $\breve{q}^{(h,1)}(\x)$ corresponds to replicating all variation in $x_1,$ for any fixed value of $x_2,$
but shifted in phase, and {\em mutatis mutandis} the analogous statements
hold for $\breve{q}^{(h,2)}(\x)$ and $\breve{q}^{(h,3)}(\x)$. If $q(\x)$ corresponds to a particular time-frequency
structure as a signal in $x_1$ for fixed values of $x_2$
then $\breve{q}^{(h,1)}(\x)$ will replicate the same structure, but shifted
in phase, in analogue with equation (\ref{phase}). We propose to use the decomposition coefficients of the three HCTs of the observed image to estimate the
decomposition of the deterministic image. An
improvement in estimation will ensue if the magnitudes of the decomposition
coefficients of the HCTs are large when the coefficient of the observed image
should be kept rather than killed. Given each coefficient replicates the
same variational structure in each separate axes this should be the case,
{\em cf} equation (\ref{phase}). For simple 1-D structures such as line segments
observed in 2-D the
energy of the image will be spread over more coefficients than strictly necessary.
By defining the additional images that should have the same marginal variational
structure as the images in each of the two axes, for moderate SNRs the estimation
should improve by using the additional components, as more often the signal
is recognized as present. Subsequent risk calculations in section \ref{riskcalcy} show that if the DWT of the quadrature components have
the same mean as the DWT of the original signal, then the risk can be reduced by a new procedure that uses the magnitudes of all four components at each
fixed $\bm{\xi},$ to threshold a given coefficient.
The authors of \cite{olhede2004} denoised a 1-D signal by defining a second
component as the HT of the observed signals, and using the DWT of this component
when thresholding the observed signal. A naive 2-D extension of this method would define
a
single extra quadrature component corresponding to $\breve{q}^{(h,3)}(\x),$ phase shifting
in both spatial directions simultaneously. We discuss the risk of this procedure
in section \ref{riskcalcy}, and it is shown to exceed that of the proposed
method, for certain scenarios. 

\section{Estimation}
\subsection{Defining Estimates}
We have argued that the quadrature components defined either
by the HCTs or RTs have the same space and spatial
frequency structure as the original image, shifted in phase. Therefore the
mean of the DWTs of the quadrature components should be the same as the DWT
of the signal. We define
a local magnitude
in terms of the DWTs of the full set of quadrature components. 

\begin{definition}[The Magnitude of a Coefficient \label{magofdef}]
{\em{We define the magnitude of a coefficient $W_{\bm{\xi}}^{(q)}$ using quadrature components denoted by
$s$ via:}}
$
M_{\bm{\xi}}^{(q,s)2}=\frac{1}{C_L^{(s)}+1}\sum_{l=0}^L W_{\bm{\xi}}^{(q,s,l)2}.
$
\end{definition}
Let for some fixed $L_1\in{\mathbb{N}}$ $B_i^j=
\left\{k:\;k=i+l,\;l=-L_1,\dots L_1\right\},$ then $M_{\bm{\xi}}^{(q,s)2}$
is a 2-D analogy to
to $S_{j,k}^2=\sum_{l \in B_k^j} W_{[j,l]}^{(q)2},$ used by Cai and Silverman \cite{Cai}[p.~132], to block threshold.
Each coefficient $W_{\bm{\xi}}^{(q)}$ will be estimated by hard
thresholding the observed coefficient depending on the value of $M_{\bm{\xi}}^{(Y,s)2}$
\begin{equation}
\widehat{W}_{\bm{\xi}}^{(q,s)}(\lambda^2)=\left\{\begin{array}{lcr}
W_{\bm{\xi}}^{(Y)} & {\mathrm{if}} & M_{\bm{\xi}}^{(Y,s)2}\ge \frac{\sigma^2}{
C_L^{(s)}+1} \lambda^{2}\\
0 & {\mathrm{if}} & M_{\bm{\xi}}^{(Y,s)2} < \frac{\sigma^2}{C_L^{(s)}+1} \lambda^{2}\\
\end{array}\right. .
\end{equation}
The notation given $\bm{\xi}$ is fixed when estimating any
set coefficient is needlessly
complicated, and we remove the reference to most of these indices. 
We define: $
\bm{W}^{(s,u)}=\left[W_{\bm{\xi}}^{(Y,s,0)},\dots,W_{\bm{\xi}}^{(Y,s,L)}\right]^T,$
$\bm{\mu}^{(s,u)}=\left[W_{\bm{\xi}}^{(q,s,0)},\dots,W_{\bm{\xi}}^{(q,s,L)}\right]^T,$
note that $\mu_1=\mu_1^{(s,u)}$ does not depend on the choice of $s$
and we denote the estimator using the $s$ indexed components by:
$\widehat{\mu}_1^{(s,u)}(\lambda^2)=\widehat{W}_{\bm{\xi}}^{(q,s)}(\lambda^2).$
%Thus
%\[\bm{W}^{(s,u)}=\bm{\mu}^{(s,u)}+\mathbf{n}^{(s,u)}.\]
%We seek to estimate $\mu_1$ based on thresholding the observed coefficients
%$W_1^{(s,u)}$
%given a full set of coefficients $\bm{W}^{(s,u)}$ are available. We shall %now
%demonstrate that for some typical signal features the magnitude of the coefficients
%is large when defined in this fashion.

\subsection{Magnitude of Typical Image Features \label{typimagefeat}}
Deterministic images
are frequently modelled as the combination of texture and contours
(see for example work by Vese and Osher \cite{Vese} modelling a signal as
a bounded variation contribution plus a texture contribution). We consider observing
an image that is an aggregation of edges and texture, where each texture
component is modelled by $t(\x)=a_t(\x_0)\cos(2\pi \f_0^T \x),$ and each
edge component is modelled by
$e(\x)=a_e(x_1)\delta(\cos(\theta)x_1+\sin(\theta)x_2-c).$ $a_e(\cdot)$ and $a_t(\cdot)$ are assumed to be slowly varying.
In general we do not expect to observe sinusoids or discontinuities
that for very slowly varying $a_e(\cdot)$ span the entire observed image,
but to be able to carry out theoretical calculations stylized image structures
must be analysed, that observed images would approximate.
We shall investigate how the magnitude of the transform coefficients of the
full set of quadrature components behave,
as our subsequent risk calculations will demonstrate that the success of
the method strongly depends on the mean of the quadrature components. 

To this purpose we define the Maximum Overlap Discrete Wavelet
Transform Coefficients (MODWT coefficients)
$\widetilde{W}_{\bm{\xi}}^{(\epsilon,s,l)}.$ These are the DWT
coefficients calculated without subsampling, and where a new normalisation
is introduced at each level $j$ to preserve energy. For a full length discussion
see Percival and Walden \cite{percival}[Ch.~4].
We denote the FT of the MODWT filter $\tilde{h}_{j,u,\kk}$
by $\tilde{H}_{j,u}(\f)=\left|\tilde{H}_{j,u}(\f)\right|e^{-2\bm{i}\pi\tilde{\varphi}_{j,u}(\f)},$
thus defining the modulus ($|\tilde{H}_{j,u}(\f)|$) and phase
($\tilde{\varphi}_{j,u}(\f)$) of $\tilde{H}_{j,u}(\f).$ For notational convenience let $\x=2^j(\kk+\bm{1})-\bm{1},$ and the region of frequency space
where $\tilde{H}_{j,u}(\f)$ is mainly supported be denoted $\Omega_{j,u}.$
The DWT coefficients of
a generic signal $q(\x)$
can be extracted from the MODWT coefficients of $q(\x),$ using the relations
$W^{(q)}_{j,u,\kk}=2^j \tilde{W}^{(q)}_{j,u,\x}$ (see for example
Percival and Walden \cite{percival}[p.~203]). 

\begin{lemma}[RT Magnitude of Local Oscillation \label{rieszosc}]
If the signal locally takes the form \\ $t(\x)=
a_t(\x_0)\cos(2\pi \f_0^T \x)$ with $\f_0=\left[f_0\cos(\phi_0)\quad f_0\sin(\phi_0)
\right]^T,$
then the magnitude of the wavelet decomposition defined by definition \ref{magofdef}
is given by:
\begin{eqnarray}
%\nonumber
%\mu_1
%&=& a_t(\x_0)\cos(2\pi (\f_0^T\x-\varphi_{j,u}(\f_0)))
%I(\f_0 \in \Omega_{j,u})+\rho_1
%\label{firstmean}\\
%\nonumber
%\mu_2^{(r,u)}
%&=& a_t(\x_0)\cos(\phi_0)\sin(2\pi (\f_0^T\x-\varphi_{j,u}(\f_0)))
%I(\f_0 \in \Omega_{j,u} )+\rho_2\\
%\mu_3^{(r,u)}
%&=& a_t(\x_0)\sin(\phi_0)\sin(2\pi (\f_0^T\x-\varphi_{j,u}(\f_0)))
%I(\f_0 \in \Omega_{j,u})+\rho_3\\
M_{\bm{\xi}}^{(t,r)2} &= &  2^{2j-1} a_t^2(\x_0) I(\f_0 \in \Omega_{j,u} )+\rho_1+O\left(\frac{1}{N}\right),
\end{eqnarray}
where $\rho_1,$ is an error term depending
on the leakage of the wavelet filters in the frequency domain. If a sufficiently
long wavelet filter is used, this term can be ignored. See Nielsen \cite{nielsen}
for more discussion on avoiding leakage.
\end{lemma}
\begin{proof}
See appendix B.
\end{proof}
\begin{lemma}[HCT Magnitude of Local Oscillation \label{hyperosc}]
If the signal locally takes the form\\ $t(\x)=
a_t(\x_0)\cos(2\pi \f_0^T \x)$ with $\f_0=\left[f_0\cos(\phi_0)\quad f_0\sin(\phi_0)
\right]^T,$
then the magnitude defined in definition \ref{magofdef} is given by:
\begin{eqnarray}
%\nonumber
%\mu_1
%&=& a_t(\x_0)\cos(2\pi (\f_0^T\x-\varphi_{j,u}(\f_0)))
%I(\f_0 \in \Omega_{j,u} )+\rho_1
%\nonumber\\
%\nonumber
%\mu_2^{(h,u)}
%&=& a_t(\x_0)\sin(2\pi (\f_0^T\x-\varphi_{j,u}(\f_0)))
%I(\f_0 \in \Omega_{j,u} )+\rho_5\\
%\nonumber
%\mu_3^{(h,u)}
%&=& a_t(\x_0)\sin(2\pi (\f_0^T\x-\varphi_{j,u}(\f_0)))
%I(\f_0 \in \Omega_{j,u})+\rho_6\\
%\nonumber
%\mu_4^{(h)}
%&=& a_t(\x_0)\cos(2\pi (\f_0^T\x-\varphi_{j,u}(\f_0)))
%I(\f_0 \in \Omega_{j,u})+\rho_7\\
M_{\bm{\xi}}^{(t,h)2} &=& 2^{2j-1} a_t^2(\x_0) I(\f_0 \in \Omega_{j,u})+\rho_2+O\left(\frac{1}{N}\right),
\label{maggynew}
\end{eqnarray}
where $\rho_2,$ is an error term depending
on the leakage of the wavelet filters in the frequency domain. 
\end{lemma}
\begin{proof}
See appendix B.
\end{proof}
\begin{lemma}[RT Magnitude of Discontinuity \label{rieszdisc}]
If the signal locally can be approximated by \\ $e(\x)=
a_e(x_1)\delta(\cos(\theta)x_1+\sin(\theta)x_2-c),$ 
then the magnitude defined in definition \ref{magofdef} is given by:
\begin{eqnarray}
\nonumber
W_{\bm{\xi}}^{(e)}&=&2^j A_e(0)\int_{\nu_{1,\min}}^{\nu_{1,\max}}
\tilde{H}_{j,u}(\cos(\theta)\nu_1,\sin(\theta)\nu_1)\;e^{2\pi \bm{i}
\nu_{1}(\cos(\theta)x_1+\sin(\theta)x_2-c)}\;d\nu_{1}+\rho_3+O\left(\frac{1}{N}\right)\\
\nonumber
U_{\bm{\xi}}^{(e)}&=&2^j A_e(0)(-\bm{i})\int_{\nu_{1,\min}}^{\nu_{1,\max}}
{\mathrm{sgn}}(\nu_1)
\tilde{H}_{j,u}(\cos(\theta)\nu_1,\sin(\theta)\nu_1)\;e^{2\pi \bm{i}
\nu_{1}(\cos(\theta)x_1+\sin(\theta)x_2-c)}\;d\nu_{1}\\
&&+\rho_4+O\left(\frac{1}{N}\right),\quad\label{eqp2}
M^{(e,r)2}_{\bm{\xi}}=\frac{1}{2}W_{\bm{\xi}}^{(e)2}+
\frac{1}{2}U_{\bm{\xi}}^{(e)2}+\rho_{5}+O\left(\frac{1}{N}\right)
\end{eqnarray}
where $\rho_{3},$ $\rho_4$ and $\rho_5$ depend on the smoothness of $a_e(\cdot),$
whilst $\nu_{1,\min}$
and $\nu_{1,\max}$ are given in appendix B.
\end{lemma}
\begin{proof}
See appendix B.
\end{proof}
\begin{lemma}[HCT Magnitude of Discontinuity \label{hyperdisc}]
If the signal locally can be approximated as a discontinuity \\
$e(\x)=
a_e(x_1)\delta(\cos(\theta)x_1+\sin(\theta)x_2-c),$ 
then the magnitude defined in definition \ref{magofdef} is given by:
\begin{eqnarray}
M^{(e,r)2}_{\bm{\xi}}=\frac{1}{2}W_{\bm{\xi}}^{(e)2}+
\frac{1}{2}U_{\bm{\xi}}^{(e)2}+\rho_{6}+O\left(\frac{1}{N}\right),
\end{eqnarray}
where $\rho_{6}$ depends on the smoothness of $a_e(\cdot),$ and the forms
of  $W_{\bm{\xi}}^{(e)}$ and $U_{\bm{\xi}}^{(e)}$ are given in lemma \ref{rieszdisc}.
%a constant envelope
%$a_e(x)=I\left(x\in\left[0,L\right)\right)/\sqrt{L}$ for example,
%will yield $\rho_{16}\equiv 0,$ as $L\rightarrow \infty.$
\end{lemma}
\begin{proof}
See appendix B.
\end{proof}
For oscillatory signals the magnitude hence aptly reflects signal presence at $\bm{\xi}.$
Equation (\ref{eqp2})
illustrates the problem experienced by an edge in a 2-D separable representation: only if $\theta=0$ or $\theta=\pi/2
m$ for $m\in {\mathbb{Z}}$ will the edge live in $\Omega_{j,2}$ or $\Omega_{j,3},$ {\em
i.e.} constant in one direction and variable in the other. Only
in this case will the representation be extremely compressed (note that the proof
needs to be adjusted for $\theta=0$).
From equation (\ref{eqp2}) we note that the compact spatial 
support of $\tilde{h}_{\bm{\xi}_{\x}}$ ensures that the
energy of $U_{\bm{\xi}}^{(e)}$ and $W_{\bm{\xi}}^{(e)}$ are mainly concentrated near $\cos(\theta)x_1+
\sin(\theta)x_2=c.$ Given we may represent $\tilde{H}_{j,u}(\cdot)$
in terms of a magnitude and a phase, the difference between $\tilde{\varphi}_{j,u}(
\cos(\theta)\nu,\sin(\theta)\nu)$
and $2\pi\nu (\cos(\theta)x_1+\sin(\theta)x_2-c)$ will determine exactly
at which spatial indices
$W_{\bm{\xi}}^{(e)}$ and $U_{\bm{\xi}}^{(e)}$ have non-negligible magnitudes
(see for example the discussion in Gopinath \cite{gopinath2005}[p.~1794]).
Thus the DWT of the quadrature components will be large near the discontinuity and they can be used to improve the estimation. Our proposed procedure will capitalise on this fact, and
it can be noticed in the reconstructions that line discontinuities are better
reconstructed (see Figure \ref{boat} (c)), and as a curved discontinuity can be approximated as the aggregation of amplitude modulated line discontinuities,
improvements in estimation can be observed for curved structures (see Figure \ref{boat} (g)).

\subsection{Distribution of Noise \& Universal Thresholds}
The distribution 
of ${\mathbf{n}}^{(s,u)}$ must be determined, to obtain a universal threshold
\cite{downie}. Let $K=N^2$ be the total number
of coefficients of
the original observed image and denote 
$M_K^{(s,\max)}=\max_{\bm{\xi}} (C_L^{(s)}+1) M_{\bm{\xi}}^{(s)2}/\sigma^2.$
Downie and Silverman \cite{downie} proposed that a universal
threshold should in general satisfy taking a value such that $\lim_{K\rightarrow \infty} P\left(
M_K^{(s,\max)}\le \lambda^2_K\right)=C,$ for some constant $0<C<1,$ and as $K$ increases the expected number of coefficients exceeding
the threshold is some small but finite non-zero value. 
We choose a slightly more conservative threshold, so that if our strategy were
to be used for $K$ independent threshold criteria, then the probability that the
maximum exceeded the threshold is 
$O\left(\frac{1}{(\log(K))^{3/2}}\right),$
rather than tending to a positive constant for the Riesz threshold. We use
for the Hypercomplex threshold a conservative version of that suggested by
Downie and Silverman. 
We cannot quite achieve
analogous results to Downie and Silverman as as the set of DWT coefficients $\mathbf{n}^{(s,u)}$
are correlated across indices $\bm{\xi},$ and adopt arguments
similar to those given by Johnstone and Silverman \cite{johnstone2} and Olhede \& Walden \cite{olhede2004},
to justify the choice of threshold.
We do not aim to determine the full covariance structure of the full set of wavelet
coefficients of the observed image and the quadrature components.
To derive the conservative threshold define ${\mathcal{M}}_K^{(s)}$ as the maximum of $K$
{\em independent} variates with the same marginal distribution as 
$\left\{ (C_L^{(s)}+1) M_{\bm{\xi}}^{(s)2}/\sigma^2\right\}.$
We find a universal threshold based on determining the distribution of ${\mathcal{M}}_K^{(s)},$ and this then constitutes
a conservative choice for $M_K^{(s,\max)}=\max_{\bm{\xi}} 
(C_L^{(s)}+1)M_{\bm{\xi}}^{(s)}/\sigma^2,$ $s=r,\;h$
as
\begin{eqnarray*}
P\left(M_K^{(s,\max)}\le \lambda^{(s)2}_K\right)&=&
P\left(\cap_{i=1}^{K} \left[ \frac{ M_{\bm{\xi}}^{(s)2}}{\sigma^2/(C_L^{(s)}+1)}\le
\lambda^{(s)2}_K\right]
\right)
\ge  \prod_{i=1}^{K}
P\left( \frac{ M_{\bm{\xi}}^{(s)2}}{\sigma^2/(C_L^{(s)}+1)}\le \lambda^{(s)2}_K\right)\\
&=& P\left({\mathcal{M}}_K^{(s)}\le \lambda^{(s)2}_K\right),
\end{eqnarray*}
by corollary 2 from Dykstra \cite{dykstra}. $C_L^{(s)}+1$ can be interpreted as the degrees of freedom
associated with $M_{\bm{\xi}}^{(s)2}.$ As shown in the subsequent section,
the DWT of the quadrature
components for any fixed value of $\bm{\xi}$ are uncorrelated, and the
$\bm{A}_i$ matrices of Dykstra are defined to take the variance
of $n_l^{(s,u)},$ into account.

\subsection{Distribution of the Magnitude}
$\var(n_{l}^{(s,u)})$ must be determined for $s=r,\;h$
to derive the approximate distribution of ${\mathbf{n}}^{(s,u)}.$ We denote
by $\overset{\cal{L}}{=}$ as equality in law \cite{ferguson}.
\begin{lemma}
[Distribution of Riesz Coefficients \label{rieszcoef}]
{\textit{The DWT coefficients of the original signal and the RTs of Gaussian
white noise are distributed as:}}
\begin{equation}
\label{covvy}
{\mathbf{n}}^{(r,u)} \overset{\cal{L}}{=}
\bm{Z}^{(r,u)}+O\left(1/N\right),\quad\bm{Z}^{(r,u)}\sim {\mathcal{N}}_3
\left(\bm{0},\sigma^2
{\mathrm{diag}}\left[
1\quad a^{(r,u)}\quad 1-a^{(r,u)}
\right]\right),\;u=1,2,3,4.
\end{equation}
where ${\mathrm{diag}}$ denotes a diagonal square matrix,
$a^{(r,u)}=\frac{1}{2},\;u=1,\;4,$
$a^{(r,2)}=\frac{1}{2}+ 2\tan^{-1}\left(\frac{1}{2}\right)+ \frac{1}{2}\tan^{-1}\left(2\right)$
and
$a^{(r,3)}=\frac{1}{2}- 2\tan^{-1}\left(\frac{1}{2}\right)- \frac{1}{2}\tan^{-1}\left(2\right).$
%\end{eqnarray}
\end{lemma}
\begin{proof}
For the proof see Appendix C.
\end{proof}
%\begin{figure*}[t]
%\centerline{
%\includegraphics[height=1.70in,width=1.70in]{idag10lf.eps}
%\includegraphics[height=1.70in,width=1.70in]{idag2lf.eps}
%\includegraphics[height=1.70in,width=1.70in]{idag4lf.eps}
%\includegraphics[height=1.70in,width=1.70in]{idag8lf.eps}}
%\caption{\label{riesz2}
%The noisy image (a), the clean image (b), the regular hard-thresholded image
%reconstruction
%(c) and the HT components hard-thresholded image
%reconstruction (d). In
%this image we have an SNR of 2.41.
%}
%\end{figure*}
\begin{lemma}
[Distribution of Riesz Magnitude \label{rieszmag}]
{\textit{The magnitude square of the DWT of the RTs of Gaussian
white noise, denoted $M^{(\epsilon,r)2}_{\bm{\xi}},$ are distributed as}}
\begin{equation}
\frac{M^{(\epsilon,r)2}_{\bm{\xi}}}{\sigma^2/2}\overset{\cal{L}}{=}
\frac{\bm{Z}^{(r,u)T}\bm{Z}^{(r,u)}}{\sigma^2}+O\left(1/N\right)
\overset{\cal{L}}{=}T_1+O\left(1/N\right),\quad
T_1 \sim \chi^2_1+a^{(r,u)}\chi^2_1+(1-a^{(r,u)})\chi^2_1,
\end{equation}
$j=1,\dots,J,$ $k_1,\;k_2=0,\dots,N_j,$
$u=1,\dots,4,$ {\textit{where if}} $u=1,\;4,$ $T_1$
{\textit{has distribution}}
\begin{equation}
\label{quatpdf}
f_{T_1}(t)=e^{-t}\frac{2}{\sqrt{\pi}}\int_{0}^{\sqrt{\frac{t}{2}}} 
e^{w^2}\;dw,
\end{equation}
{\textit{whilst if}} $u=2,\;3$ {\textit{the moment generating function of}} $T_1$ {\textit{
is readily calculable, and we may calculate the probability of obtaining
large variates
using formulae derived by \cite{grad}.}}
\end{lemma}
\begin{proof}
For the proof see Appendix C.
\end{proof}

\begin{lemma}
[Distribution of HCT Coefficients \label{hypercoef}]
{\textit{The DWT coefficients of the  HCT of Gaussian
white noise are distributed as}}
\begin{equation}
\label{covvy2}{\mathbf{n}}^{(h,u)}\overset{\cal{L}}{=} \bm{Z}^{(h,u)}+O\left(1/N\right),\quad
\bm{Z}^{(h,u)} \sim {\mathcal{N}}_4\left(\bm{0},\sigma^2
%\begin{pmatrix}
{\mathrm{diag}}\left[
1\quad 1 \quad 1 \quad 1 \right]
%1 & 0 & 0 & 0\\
%0 & 1 & 0 & 0\\
%0 & 0 & 1 & 0\\
%0 & 0 & 0 & 1
%\end{pmatrix}
\right),\;u=1,2,3,4.
\end{equation}
\end{lemma}
\begin{proof}
For the proof see Appendix C.
\end{proof}

Given the approximate joint distribution of the DWT coefficients at $\bm{\xi}$ has been determined, it trivially follows that the magnitude is distributed as 
\[\frac{M_{\bm{\xi}}^{(\epsilon,h)}}{\sigma^2/4}\overset{\cal{L}}{=}
\frac{ \bm{Z}^{(h,u)T} \bm{Z}^{(h,u)}}{\sigma^2}+O\left(1/N\right)\overset{\cal{L}}{=}
T_2+O\left(1/N\right),\quad T_2 \sim
\chi^2_1+\chi^2_1+\chi^2_1+\chi^2_1\sim
\chi^2_4.\]

\subsection{Threshold Choice}
\begin{lemma}[Riesz Conservative Threshold \label{rieszconserv}]
{\textit{Taking}} 
\begin{equation}
\label{thressy}
\lambda_K^{(r)2}(l)=2\log\left(K\right)+2C \log\left(\log\left(K\right)\right),
\end{equation}
{\textit{it follows that if}} $C>-1,$
\begin{equation}
P\left({\mathcal{M}}_K^{(r)}<\lambda_K^{(r)2}(l) \right)\rightarrow 1.
\end{equation}
\end{lemma}
\begin{proof}
For the proof see Appendix C.
\end{proof}
From \cite{downie} we may note that the RT threshold is thus
like that of a $\chi^2_1$ ($C>-1$)
however to ensure that the probability tends to
$1$ rather than a fixed constant we take $C=0,$ rather than $C=-1.$ 
Given the normalised marginal magnitudes of the HT components
are $\chi^2_4$ we may use results of \cite{downie} to note that
\[\lambda^{(h)2}_K(l)=2\log(K)+2\log(\log(K)).\]
gives an appropriate threshold.
Note that yet again, we expect this to be a conservative threshold, because
the wavelet coefficients will be correlated {\em across}
$\bm{\xi}.$ As a final step of the procedure we implement cycle-spinning
\cite[p.~429]{percival}, which is known to improve mean square error results
considerably. Finally for completeness consider implementing hard thresholding
in the usual fashion: this will be denoted by $s=c,$ and we discuss using
a single added extra component of $n_4^{(h,u)}$ when thresholding $n_1^{(h,u)}$
as a naive extension of `analytic' thresholding, denoted by taking $s=a.$

\subsection{Risk Calculations \label{riskcalcy}}
\begin{figure*}[t]
\centerline{
\includegraphics[scale=0.62]{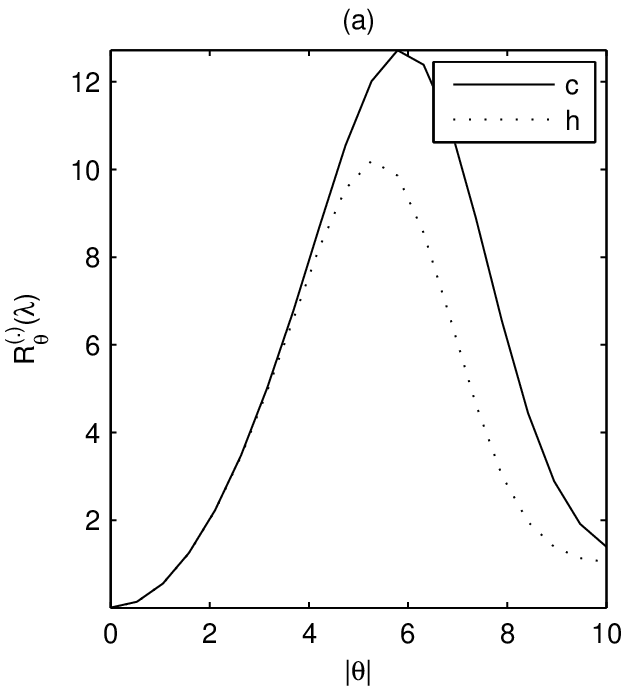}
\includegraphics[scale=0.62]{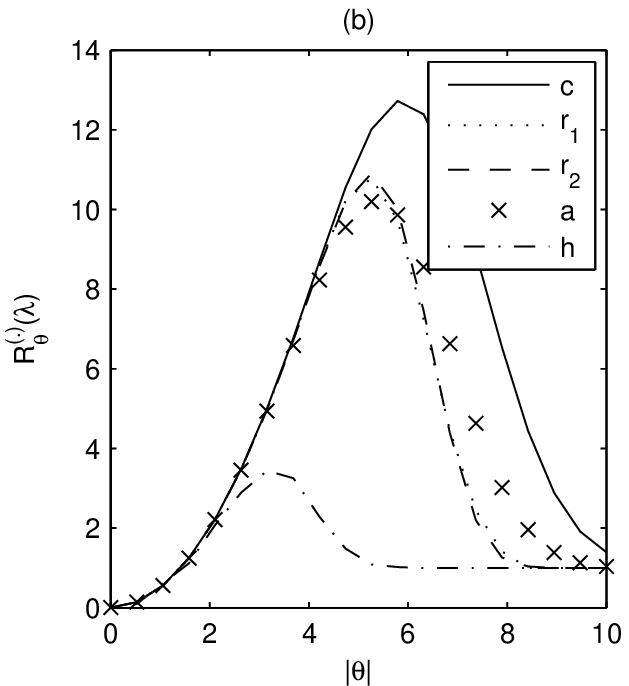}
\includegraphics[scale=0.62]{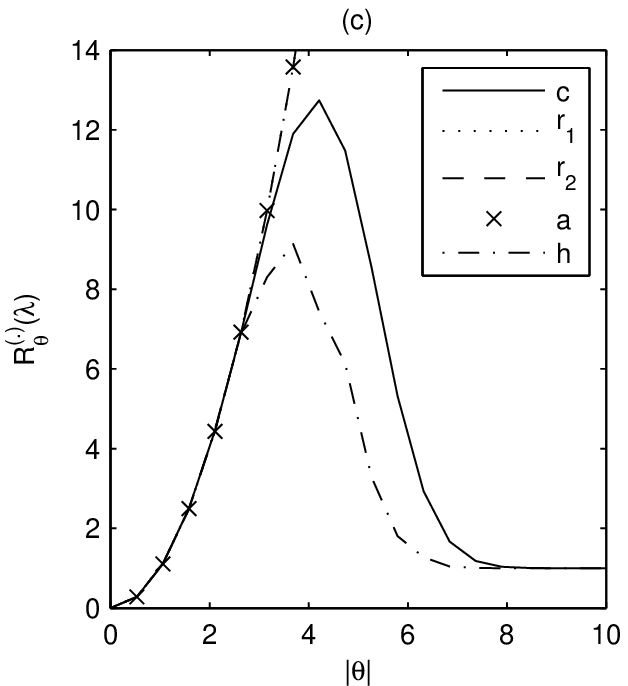}
\includegraphics[scale=0.62]{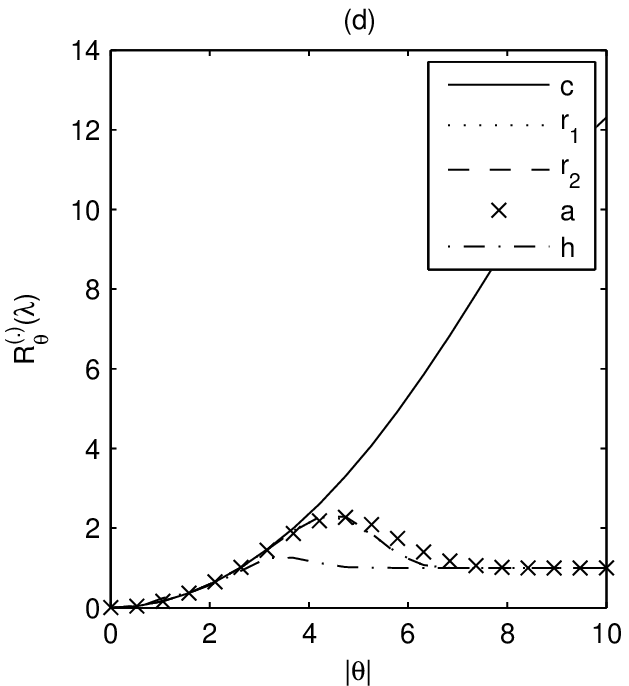}
}
\caption{\label{exx2}
The risk of hard thresholding compared to `analytic' hard thresholding (a)
where the standardised mean of the coefficient is denoted $\theta=\left|\theta\right|\cos(\phi),$
where $\phi=\pi/4.$ 
%i.e. when the signal and its HT
%have wavelet transforms of equal magnitude.
The risk associated with a thresholded coefficient using standard hard thresholding
(solid line), the analytic threshold (crosses), the Riesz thresholds (dotted line and dashed line) or the Hypercomplex threshold
(dash-dotted line). In plot (b) $\theta_1=\sqrt{\theta_2^2+\theta_3^2}=\left|\theta\right|/\sqrt{2}$ for
the Riesz threshold whilst $\theta_1=\theta_2=\theta_3=\theta_4=\left|\theta\right|/\sqrt{2},$ for the Hypercomplex
threshold. In plot (c) $\theta_2=\theta_3=0$ for the Riesz threshold whilst
$\theta_1=\theta_2$ and
$\theta_3=\theta_4=0$ for the Hypercomplex threshold. In plot (d) $\theta_1=
\left|\theta\right|\cos(3\pi/8)$ 
$\theta_2=\theta_3=
\left|\theta\right|\sin(3\pi/8)/\sqrt{2}$ for the Riesz threshold whilst
$\theta_2=\theta_1=\left|\theta\right|\cos(3\pi/8)$ and
$\theta_3=\theta_4=\left|\theta\right|\sin(3\pi/8)$ for the Hypercomplex threshold.  For the `analytic' procedure we use $\theta_4$ as the second
component.
%In plot (d) $\theta_1=
%\left|\theta\right|\cos(\pi/3)$ 
%$\theta_2=\theta_3=
%\left|\theta\right|\sin(\pi/3)/\sqrt{2}$  for the Riesz threshold whilst
%$\theta_2=\theta_1=\left|\theta\right|\cos(\pi/3)$ and
%$\theta_3=\theta_4=\left|\theta\right|\sin(\pi/3)$ for the Hypercomplex %threshold.
}
\end{figure*}

To compare the theoretical performance of the threshold estimators proposed
in this paper, we calculate the standardized mean
square risk
at any fixed value of $\bm{\xi}.$ We define the standardized risk using any
threshold procedure denote by $s$ for $\theta_l=\mu_l^{(s,u)}/\sigma$ by
\begin{equation}
R_{\theta}^{(s)}(\lambda)=\sigma^{-2}
E\left[\left(\widehat{\mu}_1^{(s,u)}(\lambda^2)-\mu_1\right)^2\right].
\end{equation}
If $s=r$ then we denote by $r_1$ and $r_2$ the two different cases that may
occur at a given $\bm{\xi}$ when $u=1,4$ or $u=2,3$ --
the risk will be different in these two cases. This will not happen for
$s=c$ or $s=h.$
For completeness we here also provide the risk of the `analytic' denoising,
as this was not done in Olhede \& Walden \cite{olhede2004} and corresponds
to a special case of the risks determined by Fryzlewicz \cite{fryz}.
\begin{theorem}
[The Risk of a Thresholded Coefficient \label{apy1}]
{\textit{The standardized risk of an individual coefficient is using threshold strategy
$s=c,\;a,\;r,\;h$ with $\theta_i=\mu_i^{(s,u)}/\sigma$ and 
$R_j(\lambda)=\sum_l(w_l+\theta_l)^2< \lambda^2$ given by:}}
\begin{eqnarray}
\nonumber
R_{\theta}^{(c)}(\lambda)&=&1+\int_{(w+\theta_1)^2< \lambda^2}\left[\theta_1^2-w^2\right]
\phi(w)\;dw,\;
R_{\theta}^{(a)}(\lambda)=1+\int_{R_j(\lambda)}\left[\theta_1^2-w^2\right]
\phi(w_1)\phi(w_2)\;d^2\bm{w}\\
\nonumber
R_{\theta}^{(r)}(\lambda)&=&1+\frac{1}{\sqrt{a^{(r,u)}(1-a^{(r,u)})}}
\int_{R_j(\lambda)}\left[\theta_1^2-w^2\right]
\phi(w_1)\phi(w_2/\sqrt{a^{(r,u)}})\phi(w_3/\sqrt{1-a^{(r,u)}})\;d^3\bm{w}\\
R_{\theta}^{(h)}(\lambda)&=&1+\int_{R_j(\lambda)}\left[\theta_1^2-w^2\right]
\phi(w_1)\phi(w_2)\phi(w_3)\phi(w_4)\;d^4\bm{w}.
\end{eqnarray}
\end{theorem}
\begin{proof}
The risk of an individual coefficient using standard hard thresholding has
been noted by Marron {\em et al.} \cite{Marron},
whilst the risk of the hyperanalytic thresholds are derived in appendix D,
and $s=a$ is a special case of the bi-variate thresholding investigated by
Fryzlewicz \cite{fryz}.
\end{proof}
For some examples of signal/noise distributions, the individual risk of a
given coefficient is plotted in Figure
\ref{exx2} for the four estimation procedures using the
universal threshold. Figure \ref{exx2} (a) shows the reduced risk of `analytic'
thresholding compared to regular thresholding in 1-D when the means of the wavelet
coefficient of the signal and of the HT of the signal are equal, and this
then
provides theoretical justification for the `analytic' denoising procedure.
Figures
\ref{exx2} (b), (c) and (d) show the risk associated with thresholding at
any index $\bm{\xi}$ using either the usual hard thresholding ($c$), `analytic'
denoising ($a$), Riesz denoising when $u=1,4$ ($r_1$) and $u=2,3$ ($r_2$)
or Hypercomplex denoising ($h$). The risk is calculated with $K=256^2$ and
using the universal threshold. If the mean of the DWT of the quadrature components
is of similar magnitude to the DWT of the signal then the risk is reduced. The greatest
weakness of the proposed methods is if the means of the quadrature components are completely
disparate from that of the original signal, as may be noted from figure \ref{exx2}
(c). The results of section \ref{typimagefeat} indicate that this will not
be the case for typical image features. The norm of a signal and its
HT are identical, and given the results of section \ref{typimagefeat} the means are unlikely to be consistently
mismatched. Finally if there is no signal present we observe the following
result.
\begin{corollary}
[The Risk of a Thresholded Coefficient when there is no signal \label{risk0}]
{\textit{The risk of an individual coefficient is using threshold strategy
$s=c,\;a,\;r,\;h$ with $\theta_{l+1}=0$ for $l=0,\dots,L$ given by:}}
\begin{eqnarray}
\nonumber
R_{0}^{(c)}(\lambda)&=&1-\gamma\left(\frac{1}{2},\frac{1}{2}\lambda^2\right),\;
R_{0}^{(a)}(\lambda)=e^{-\frac{1}{2}\lambda^2}\left(1+\frac{1}{2}\lambda^2\right),\;
R_{0}^{(h)}(\lambda)=e^{-\frac{1}{2}\lambda^2}\left(1+\frac{1}{2}\lambda^2+
\frac{1}{8}\lambda^4\right)\\
R_{0}^{(r)}(\lambda)&=&1+\gamma\left(\frac{1}{2},\frac{1}{2}\lambda^2\right)-2
\gamma\left(\frac{1}{2},\lambda^2\right)+4\lambda^2\left[\Phi(\sqrt{2}\lambda)-
\Phi(\lambda)\right],\;u=1,\;4.
\end{eqnarray}
\end{corollary}
\begin{proof}
See appendix D. We denote by $\gamma(a,x)=\frac{1}{\Gamma(a)}\int_{0}^x s^{a-1} e^{-s}\;ds.$
\end{proof}
As the representation of the image will be sparse the risk if no signal is
present is important.
From the corollary, and the asymptotic forms in $\lambda$ given in appendix D, we may note that the risk at the universal threshold when there
is no signal present is of the
same order for $R_{0}^{(c)}(\sqrt{2\log(K)})$ and 
$R_{0}^{(r_1)}(\sqrt{2\log(K)})$
even if the coefficient differs in favour of $R_{0}^{(c)}(\sqrt{2\log(K)})$
whilst $R_{0}^{(h)}(\sqrt{2\log(K\log(K))})$ and \\ $R_{0}^{(a)}(\sqrt{2\log(K)})$
correspond to different orders. The thresholds were introduced
to improve the estimation of signals that were slightly more spread across coefficients
than strictly necessary, but the risk for any coefficient when no
signal is present is of similar enough nature to make
the difference in estimation negligible (i.e. $O\left(\frac{\log(K)}{K}\right)$ rather
than $O\left(\frac{1}{K\log(K)}\right).$ The examples will substantiate this
claim.

\section{Examples}
To examine the properties of the proposed methods, we have implemented simulation
studies on images that can be retrieved at \verb=http://sipi.usc.edu/database/=
(Tiffany and Boat), whilst (Lenna and MRIScan) are downloaded from 
\verb=http://www-stat.stanford.edu/~wavelab/=. We used LA wavelets length 8.
To compare our results, similarly to \cite{starck},
we also implemented usual hard thresholding and the
Wavelet-domain Hidden Markov Model (HMM method) proposed by the Rice group
\cite{crouse}, where the
software is available at
\verb=http://www-dsp.rice.edu/software/=, denoted by $s=hmm.$ We used the
code ${\mathrm{hdenoise.m}}$ with default settings, and ${\mathrm{daubcqf(8,'min')}}.$
We implemented the method at several
Signal-to-Noise Ratio (SNRs) of 2 (very noisy), 4 and 8 (quite clean), with a set of images, {\em i.e.} Lenna ($512
\times 512$ version), Boat ($512 \times 512$), MRIScan ($256 \times 256$)
and the second channel of the colour image of Tiffany ($512 \times 512$). The SNR is (as usual) given
by $SNR^2=\frac{1}{N^2 \sigma^2}\sum\sum q_{\x}^2.$  Table \ref{tab:mse}
shows the result over repeated simulations. Reduced Mean Square Errors
MSEs, and increased Peak Signal to Noise Ratios PSNRs (using the definition
of \cite{barber2}), are observed when using the proposed method with
either hyperanalytic threshold criterion,
and the reduction
in MSE is of a respectable magnitude compared to variation across replications
as the estimated standard deviations in the MSEs are usually considerably smaller. Overall the hypercomplex thresholding
procedure is outperforming the Riesz thresholding as well as the other methods, apart from the boat image
at high SNRs where the HMM does better. The Hypercomplex method is expected
to outperform the Riesz method
from the risk calculations, but not perhaps from our discussion in sections
\ref{riessy} and \ref{hyperc}. The Riesz transform may appear more useful
as it determines the prevalent direction
from the image, whilst the Hypercomplex transform simply decreases the risk
in estimation by considering variation associated with the same time-frequency
(i.e. 1-D)
behaviour in both axes separately, with the second variable treated as fixed. However, whilst the Riesz transform is
suitable to use on locally unidirectional structure as discussed in section
\ref{riessy}, the hypercomplex transform
treats variation in both axes, and images quite frequently have multi-directional
variations present even locally. Some additional analysis of images has
also been implemented in \cite{olhede2006}.
Consider two
cuts from reconstructions to further elucidate on these results: see Figure \ref{boat} (a)-(h). We show the hypercomplex
reconstructions only as the Riesz and hypercomplex reconstructions are similar. Clearly both \ref{boat} (c) is more connected than (b),
as is (g) to (f), whilst (d) and (h) has much remaining noise to preserve
more detail. The proposed method performs quite
well.
In addition to the SNR's chosen for the full range of images, we implemented
the procedure at the SNR chosen by Starck {\em et al.} \cite{starck}, namely
adding Gaussian noise with a standard deviation of 20 to the raw Lenna image, or using a SNR of 5.58. The PSNR we observed in the noisy image (21.58) is
less than theirs
(22.13) but that is to be expected in noisy replications. We found
for the methods tested in this paper that averaged over 100 replications $s=c$ (29.22), $s=r$ (30.12), $s=h$ (30.93) and
$s=hmm$ (30.48), where they obtained for $s=c$ (28.35) and $s=hmm$ (30.80).
Starck obtained PSNRs between 29.99 to 31.95 by using local ridgelets
and curvelets. Clearly
the hypercomplex denoising performs on par with the algorithms suggested,
and in addition the proposed procedure is both cheap to implement and extremely
simple to code. 

\begin{table*}[ht]
\centering{
\caption{The average results over 50 runs. The symmlet wavelets (or LA wavelets
were used) and $J=3.$
\label{tab:mse}}
\begin{tabular}[b]
{|l|c|c|c|c|}
\hline
\hline
Example (SNR)                    & Boat (2/4/8)    & Lena (2/4/8)     
                                 & Tiffany (2/4/8)   & MRIScan (2/4/8)\\ \hline\hline
average MSE (n)$\times 10^{-8}$  & 95.4/23.8/6.0    &95.4/23.8/6.0
                                 & 95.4/23.8/6.0    &381.5/95.4/23.8\\
sd MSE (n) $\times 10^{-10}$     & 26.3/6.6/1.6     &26.3/6.6/1.6
                                 & 26.3/6.6/1.6     &187.9/47.0/11.7\\
PSNR (n)                         & 11.36/17.38/23.40&10.67/16.69/22.71
                                 & 7.65/13.67/19.69 &16.96/22.98/29.00\\
\hline\hline
average MSE (c)$\times 10^{-8}$  & 8.0/4.4/2.2      &5.16/2.6/1.1
                                 & 3.8/2.2/1.0      &56.9/24.5/9.1\\
sd MSE (c) $\times 10^{-10}$     & 6.7/3.2/1.2      &7.3/2.3/0.8
                                 & 4.9/2.5/0.8      &110.0/43.4/15.3\\
PSNR (c)                         &22.12/24.68/27.70 &23.34/26.34/29.88
                                 &21.64/24.02/27.39 &25.22/28.88/33.18\\
\hline\hline
average MSE (r)$\times 10^{-8}$  &7.1/3.8/1.8       &4.5/2.2/0.93
                                 &3.4/1.9/0.8       &46.9/19.7/7.3\\
sd MSE (r) $\times 10^{-10}$     &6.2/2.7/1.0       &5.6/2.1/0.7
                                 &3.8/2.4/0.7       &89.2/36.9/11.6\\
PSNR (r)                         &22.63/25.39/28.49 &23.96/27.09/30.77
                                 &22.09/24.73/28.00 &26.06/29.82/34.17\\
\hline\hline
average MSE (h)$\times 10^{-8}$  &{\bf 6.3}/3.2/1.5 &{\bf 3.9}/{\bf 1.9}/{\bf
0.77}
                                 &{\bf 3.1}/{\bf 1.7}/{\bf 0.8}&{\bf 40.4}/{\bf
                                 16.6}/{\bf 6.1}\\
sd MSE (h) $\times 10^{-10}$     &6.6/2.2/0.9       &4.6/2.0/0.6
                                 &3.3/2.3/0.6       &79.1/31.0/8.7\\
PSNR (h)                         &{\bf 23.14}/ 26.04/29.2  & {\bf 24.56}/{\bf
                                 27.74}/{\bf 31.58}
                                 &{\bf 22.48}/{\bf 25.06}/{\bf 28.55}&{\bf
                                 26.71}/{\bf 30.57}/{\bf 34.92}\\
\hline\hline
average MSE (hmm)$\times 10^{-8}$&6.7/{\bf 3.0}/{\bf 1.4}&6.3/2.1/ 0.85
                                 &5.4/1.8/{\bf 0.8} &48.9/20.9/7.5\\
sd MSE (hmm) $\times 10^{-10}$   &11.3/1.8/0.7      &16.0/2.3/0.5
                                 &16.4/1.2/2.7      &82.1/56.9/7.4\\
PSNR (hmm)                       &22.90/{\bf 26.38}/{\bf 29.62} &22.48/27.31/31.19
                                 &20.12/24.78/28.26 &25.88/29.57/34.01\\
\hline
\end{tabular}}
\end{table*}

\section{Conclusions}
This paper has proposed a new thresholding strategy for estimating decomposition
coefficients,
and has in particular implemented the strategy with the discrete separable
DWT. We
have determined the stochastic properties of the decomposition of the noise, and both the deterministic (for some stylized image features) and
stochastic properties
of the suggested new thresholding criterion. We established universal thresholds.
We calculated the risk theoretically,
and for some specific choices of the mean provided plots of the risk showing
that the proposed methods outperform standard denoising theoretically. 
We implemented the procedure on several examples at several SNRs, comparing the methods with the Hidden-Markov-Model
used by the Rice group as well as standard hard thresholding, and found that the proposed algorithms offered
improvements in most cases. Given the simplicity in implementation, and visually
pleasing reconstructions, hyperanalytic denoising methods offer a computationally
cheap improvement to existing methodology, as well as offers insight into
2-D variational structure.

\section*{Acknowledgments}
SO would like to express her thanks to the anonymous referees for
the many helpful suggestions that substantially improved the paper, as well
as her understanding of the topic. SO would also like to thank Professor Andrew Walden for introducing her to this
research area, and gratefully acknowledges financial support from EPSRC (UK).
SO gratefully acknowledges usage of WaveLab routines and the data sets.

%\bibliography{hyperbib2}

\appendix 
\section*{A: Digital Implementation}
For future reference the Discrete Fourier Transform (DFT) and its inverse
are given with $\bm{\Omega}=\left[-\frac{1}{2\Delta x},\frac{1}{2\Delta x}\right]^2$
by:
\begin{equation}
\label{fftfor}
Q_D\left(f_1,f_2\right)=\Delta x^2
\sum_{x_1=0}^{N-1} \sum_{x_2=0}^{N-1} q(x_1\Delta x,x_2\Delta x) e^{-2i\pi
\f^T\x \Delta x},\quad
q_{\x}=\int\int_{\bm{\Omega}}
Q_D(\f) e^{2i\pi
\f^T \x}\;d^2\f.
\end{equation}
Note that 
the value of $Q_D\left(f_1,f_2\right)$
at frequencies $f_1\in \left\{\frac{N/2}{N
\Delta x},\; \frac{N/2+1}{N \Delta x},\;\frac{N/2+2}{N \Delta x},\;\dots,
\;\frac{N-1}{N \Delta x}\right\}$ is equivalent to 
the value of $Q_D\left(f_1,f_2\right)$
at frequencies $f_1\in \left\{-\frac{N/2}{N
\Delta x},\; -\frac{N/2-1}{N \Delta x},\;-\frac{N/2-2}{N \Delta x},\;\dots,
\;-\frac{1}{N \Delta x}\right\},$ and the equivalent statement {\it mutatis
mutandis} hold for $f_2.$ Let $Q_{D,\mathbf{u}}=Q_D\left(\frac{u_1}{N\Delta x},\frac{u_2}{N\Delta x} \right).$
\begin{eqnarray*}
\breve{q}^{(s,l)}_{\x}&=&
\int
\int_{\bm{\Omega}}
Q_D(\f) V^{(s,l)}_D(\f) e^{2i\pi
\f^T\x \Delta x}\;d^2 \f
%&=& \frac{1}{\Delta x^2 N^2}
%\sum_{u_1=-(N/2-1)}^{N/2}\sum_{u_2=-(N/2-1)}^{N/2}
% Q_{u_1,u_2} V^{(s,l)}_D\left(\frac{u_1}{N\Delta x},
%\frac{u_2}{N\Delta x}\right)\\
%&&e^{2i\pi
%{\mathbf{u}}^T\x}+O\left(\frac{1}{N}\right)\\
=\frac{1}{\Delta x^2 N^2}
\sum_{u_1=0}^{N-1}\sum_{u_2=0}^{N-1}
 Q_{D,\mathbf{u}} V^{(s,l)}_{D,\mathbf{u}} e^{2i\pi
{\mathbf{u}}^T\x}+O\left(\frac{1}{N}\right),
\end{eqnarray*}
$s=r,\;h,$ and $l=1,\;\dots,\;L,$
where we define:
\[V^{(r,1)}_{D,\mathbf{u}}=\left\{\begin{array}{ccc}
0 &{\mathrm{if}} & u_1=0\\
\frac{u_1/(N\Delta x)}{i\sqrt{(u_1/(N\Delta x))^2+
(l_2/(N\Delta x))^2}}&{\mathrm{if}} & u_1=1,\dots,N/2,\;u_2=1,\dots,N/2\\
\frac{(N-u_1)/(N\Delta x)}
{i\sqrt{((N-u_1)/N\Delta x)^2+
(u_2/N\Delta x)^2}}&{\mathrm{if}} & u_1=N/2+1,\dots,N-1,\;u_2=1,\dots,N/2\\
\frac{u_1/(N\Delta x)}{i\sqrt{(u_1/(N\Delta x))^2+
((N-u_2)/(N\Delta x))^2}}&{\mathrm{if}} & u_1=1,\dots,N/2,\;u_2=N/2+1,\dots,N-1\\
\frac{(N-u_1)/(N\Delta x)}
{i\sqrt{((N-u_1)/N\Delta x)^2+
((N-u_2)/N\Delta x)^2}}&{\mathrm{if}} & u_1=N/2+1,\dots,N-1,\;u_2=1,\dots,N/2
\end{array}\right. .\]
We define $V^{(r,2)}_{D,\mathbf{u}}$ in an analogous fashion. The
digital definition of the filters corresponding to the three Hypercomplex
components are defined via:
\[V^{(h,l)}_{D,\mathbf{u}}=\left\{\begin{array}{ccc}
0 &{\mathrm{if}} & u_l=0,N/2\\
-i&{\mathrm{if}} & u_l=1,\dots,N/2-1\\
i&{\mathrm{if}} & u_l=N/2+1,\dots,N-1
\end{array}\right.,\;l=1,2 ,\]
and $V^{(h,3)}_{D,\mathbf{u}}=V^{(h,1)}_{D,\mathbf{u}}
V^{(h,2)}_{D,\mathbf{u}}.$ Implementing the discrete HCT introduces an error term of $O\left(1/N\right).$

\section*{B: DWT of Typical Deterministic Image Features}
We use
period boundary treatment (see Mallat \cite{mallat}[p.~282--292]) when implementing
the MODWT, so that with $\bm{\xi}_{\x}=\left[j,u,\x\right]^T:$
\begin{equation}
\tilde{W}_{\bm{\xi}_{\x}}^{(q)}=
\int\int_{\bm{\Omega}} \tilde{H}_{j,u}(\f)Q_D(\f)e^{2\pi \bm{i} \f^2\x}\;d^2\f.
\label{wavt}
\end{equation}

\subsection*{Proof of Lemma \ref{rieszosc}}
By direct calculation using equation (\ref{wavt}):
\begin{eqnarray*}
\nonumber
\mu_1&=&2^{j} \int \int_{\bm{\Omega}} \frac{a_t(\x_0)}{2}\left(\delta(\f_0-\f)+
\delta(\f_0+\f)\right)\tilde{H}_{j,u}(\f) e^{2\bm{i}\pi \f^T\x}\;d^2\f+O\left(\frac{1}{N}\right)\\
&=&2^j a_t(\x_0)\left|\tilde{H}_{j,u}(\f_0)\right|
\cos(2\pi (\f_0^T\x-\varphi_{j,u}(\f_0)))+O\left(\frac{1}{N}\right)
%\overset{(1)}{=} a_t(\x_0)\cos(2\pi (\f_0^T\x-\varphi_{j,u}(\f_0)))
%I(\f_0 \in \Omega_{j,u})+\rho_2
\label{firstmean2}\\
\nonumber
\mu_2^{(r,u)}&=&2^{j} \int\int_{\bm{\Omega}}
 \frac{a(\x_0)}{2}(-\bm{i}){\mathrm{sgn}}(f_1)\frac{f_1}{f}\left(\delta(\f_0-\f)+
\delta(\f_0+\f)\right)\tilde{H}_{j,u}(\f) e^{2\bm{i}\pi \f^T\x}\;d^2\f\\
&&+O\left(\frac{1}{N}\right)=2^j a_t(\x_0)\left|\tilde{H}_{j,u}(\f_0)\right|
\cos(\phi_0) \sin(2\pi (\f_0^T\x-\varphi_{j,u}(\f_0)))\nonumber +O\left(\frac{1}{N}\right),\end{eqnarray*}
%&\overset{(2)}{=}& a_t(\x_0)\cos(\phi_0)\sin(2\pi (\f_0^T\x-\varphi_{j,u}(\f_0)))
%I(\f_0 \in \Omega_{j,u})+\rho_2\\
and similarly 
$\mu_3^{(r,u)}=2^j a_t(\x_0)\left|\tilde{H}_{j,u}(\f_0)\right|
\sin(\phi_0) \sin(2\pi (\f_0^T\x-\varphi_{j,u}(\f_0)))+O\left(\frac{1}{N}\right),$ %\nonumber% \\
%&\overset{(3)}{=}& a_t(\x_0)\sin(\phi_0)\sin(2\pi (\f_0^T\x-\varphi_{j,u}(\f_0)))
%I(\f_0 \in \Omega_{j,u})+\rho_3\\
%M_{\bm{\xi}}^{(q,r)2} &=&  a_t^2(\x_0) I(\f_0\in \Omega_{j,u} )+\rho_4.
%\end{eqnarray*}
whilst $\rho_1$ is introduced when $\left|\tilde{H}_{j,u}(\f_0)\right|$
is approximated by an exact band-pass structure.
\subsection*{Proof of Lemma \ref{hyperosc}}
By direct calculation, using equation (\ref{wavt}), it follows that:
\begin{eqnarray*}
%\nonumber
%\mu_1
%&=& a_t(\x_0)\cos(2\pi (\f_0^T\x-\varphi_{j,u}(\f_0)))
%I(\f_0 \in \Omega_{j,u})+\rho_1
%\nonumber\\
\mu_2^{(h,u)}&=&2^{j} \int \int_{\bm{\Omega}} \frac{a(\x_0)}{2}(-\bm{i})
{\mathrm{sgn}}(f_1)\left(\delta(\f_0-\f)+
\delta(\f_0+\f)\right)\tilde{H}_{j,u}(\f) e^{2\bm{i}\pi \f^T\x}\;d^2\f\\
&&+O\left(\frac{1}{N}\right)=2^j a_t(\x_0)\left|\tilde{H}_{j,u}(\f_0)\right|
\sin(2\pi (\f_0^T\x-\varphi_{j,u}(\f_0)))\nonumber+O\left(\frac{1}{N}\right) .\end{eqnarray*}
Similarly
%&\overset{(1)}{=}& a_t(\x_0)\sin(2\pi (\f_0^T\kk-\varphi_{j,u}(\f_0)))
%I(\f_0 \in \Omega_{j,u})+\rho_5\\
$\mu_3^{(h,u)}= 2^j a_t(\x_0)\left|\tilde{H}_{j,u}(\f_0)\right|
\sin(2\pi (\f_0^T\x-\varphi_{j,u}(\f_0)))+O\left(\frac{1}{N}\right)$ %\nonumber %\\
%&\overset{(2)}{=}& a_t(\x_0)\sin(2\pi (\f_0^T\x-\varphi_{j,u}(\f_0)))
%I(\f_0 \in \Omega_{j,u})+\rho_6\\
and \\ $\mu_4^{(h,u)}=2^j a_t(\x_0)\left|\tilde{H}_{j,u}(\f_0)\right|
\cos(2\pi (\f_0^T\x-\varphi_{j,u}(\f_0)))+O\left(\frac{1}{N}\right).$ %\nonumber% %\\
%&\overset{(3)}{=}& a_t(\x_0)\cos(2\pi (\f_0^T\x-\varphi_{j,u}(\f_0)))
%I(\f_0 \in \Omega_{j,u})+\rho_7\\
%M_{\bm{\xi}}^{(q,h)2} &=& 2 a_t^2(\x_0) I(\f_0 \in )+\rho_8.
%\end{eqnarray*}

\subsection*{Proof of Lemma \ref{rieszdisc}}
Assume for simplicity of exposition, $0< \theta\le
\pi/2,$ but with the necessary notational changes no such
restriction needs to be made. Define the rotation matrix by
$\bm{r}_{\theta}=
\left[
\cos(\theta) \; -\sin(\theta),
\sin(\theta) \; \cos(\theta)
\right],$
and the change of variable given by:
$\bm{\nu}
=\f(-\theta)=\bm{r}_{-\theta}\f.$ 
We assume that $A_e(f)$ decays for large frequencies, and consider $A_e(f)= A_e(0)\delta(f-0).$ For example if
$a_e(x)=I\left(x\in\left[0,L\right)\right)/\sqrt{L}$ then $A_e(f)=\frac{e^{-i\pi
f L}\sin(\pi f L)}{\pi f \sqrt{L}},$ that
as $L \rightarrow \infty$ $A_e(f)$ will concentrate to $f=0.$
$
E(\f)=\frac{A_e\left(-f_2(\theta)/\sin(-\theta)\right)}{\left|\sin(\theta)\right|}
e^{-2\pi \bm{i}c f_2/\sin(\theta)},
$ and with $\bm{\Omega}(\theta)$ the rotated by $\theta$ version of $\bm{\Omega},$
with $\nu_{1,\min}=\min_{\nu} \nu=\sec(\theta) f_1 I\left(\left|\tan (\theta)
f_1\right|<\frac{1}{2\Delta x}\right)$ $I\left(\left|f_1\right|<\frac{1}{2\Delta x}\right)$ and $\nu_{1,\max}=\max_{\nu} \nu=\sec(\theta) f_1 I\left(\left|\tan (\theta)
f_1\right|<\frac{1}{2\Delta x}\right)I\left(\left|f_1\right|<\frac{1}{2\Delta x}\right),$
then
\begin{eqnarray}
\nonumber
\mu_1&=&2^j \int\int_{\bm{\Omega}}\tilde{H}_{j,u}(\f)\frac{A_e\left(-f_2(-\theta)/\sin(\theta)\right)}
{\left|\sin(\theta)\right|}e^{-2\pi \bm{i}c f_2/\sin(\theta)}\;e^{2\pi \bm{i}
\f^T\x}\;d^2\f+O\left(\frac{1}{N}\right)\\
\nonumber
&=&\int \int_{\bm{\Omega}(\theta)}\tilde{H}_{j,u}(\bm{r}_{-\theta}
\bm{\nu})\frac{2^j A_e\left(-\nu_{2}/\sin(\theta)\right)}
{\left|\sin(\theta)\right|}e^{-2\pi  \bm{i}((\sin(\theta)\nu_{1}+\cos(\theta)\nu_{2})
c/\sin(\theta)-\bm{\nu}^T\bm{r}_{-\theta}\x)}\;d^2\bm{\nu}+O\left(\frac{1}{N}\right)\\
&\overset{(1)}{=}& 2^j A_e(0)\int_{\nu_{1,\min}}^{\nu_{1,\max}}
\tilde{H}_{j,u}(\cos(\theta)\nu_1,\sin(\theta)\nu_1)\;e^{2\pi \bm{i}
\nu_{1}(\cos(\theta)x_1+\sin(\theta)x_2-c)}\;d\nu_{1}+\rho_3+O\left(\frac{1}{N}\right).
\nonumber
\end{eqnarray}
The approximation in (1) relies on $A_e(\cdot)$ taking the form of a $\delta$
distribution contribution, i.e. $a_e(\cdot)$ constant over a large spatial
domain, but a slowly
varying $a_e(\cdot)$ will approximately yield the same result. 
Also we may find approximate descriptions for decomposition of the RTs, namely
with $f=\sqrt{f_1^2+f_2^2}$ and $\nu=\sqrt{\nu_1^2+\nu_2^2}:$
\begin{eqnarray}
\nonumber
\mu_2^{(r,u)}&=&2^j \int\int_{\bm{\Omega}}
\tilde{H}_{j,u}(\f)(-\bm{i})f_1/f\frac{A_e\left(-f_2(-\theta)/\sin(\theta)\right)}
{\left|\sin(\theta)\right|}e^{-2\pi \bm{i}c f_2/\sin(\theta)}\;e^{2\pi \bm{i}
\f^T\x}\;d^2\f\\
&&+O\left(\frac{1}{N}\right)
\nonumber
%&=&2^j \int \int\tilde{H}_{j,u}(\bm{r}_{\theta}\bm{\nu})
%\frac{\cos(\theta)\nu_{1}-\sin(\theta)\nu_2)}{\nu}
%\frac{A_2\left(-\nu_{2}/\sin(\theta)\right)}
%{\left|\sin(\theta)\right|}e^{-2\pi c \bm{i}(\sin(\theta)\nu_{1}+\cos(\theta)\nu_{2})
%/\sin(\theta)}\;e^{2\pi \bm{i}
%\bm{\nu}^T\bm{r}_{-\theta}\x}\;d^2\bm{\nu}\\
\overset{(1)}{=} A_e(0)\nonumber 2^j(-\bm{i})
\cos(\theta) \int_{\nu_{1,\min}}^{\nu_{1,\max}}\tilde{H}_{j,u}(\cos(\theta)\nu_1,\sin(\theta)\nu_1)
\;e^{2\pi \bm{i}
\nu_{1}(\cos(\theta)x_1+\sin(\theta)x_2-c)}\\
&&\times 
{\mathrm{sgn}}(\nu_1)\;d\nu_{1}+\rho_{4}^{\prime}
+O\left(\frac{1}{N}\right)=  
\cos(\theta)U_{\bm{\xi}}^{(e)}+\rho_{4}^{\prime}+O\left(\frac{1}{N}\right).
\end{eqnarray}
The approximation in (1) relies on the envelope being constant in the spatial
domain -- a slowly
varying envelope will thus only approximately yield the same value, and this
introduces an error term $\rho_4^{\prime}.$
Similarly it transpires that
$
\mu_3^{(r,u)}=   \sin(\theta)U_{\bm{\xi}}^{(e)}+\rho_{7}+O\left(\frac{1}{N}\right),
$
and thus the result follows, with a new error term $\rho_5.$

\subsection*{Proof of Lemma \ref{hyperdisc}}
Also we may find approximate descriptions for the Hypercomplex components, namely:
\begin{eqnarray}
\nonumber
\mu_2^{(h,u)}&=&2^j \int\int_{\bm{\Omega}}
\tilde{H}_{j,u}(\f)(-\bm{i}){\mathrm{sgn}}\left(f_1\right)
\frac{A_e\left(-f_2(\theta)/\sin(-\theta)\right)}
{\left|\sin(\theta)\right|}e^{-2\pi \bm{i}c f_2/\sin(\theta)}\;e^{2\pi \bm{i}
\f^T\x}\;d^2\f\\
&&+O\left(\frac{1}{N}\right)
\nonumber
%&=&2^j \int \int\tilde{H}_{j,u}(\bm{r}_{\theta}\bm{\nu})
%(-\bm{i}){\mathrm{sgn}}\left(\cos(\theta)\nu_{1}- \sin(\theta) \nu_{2}\right)
%\frac{A_e\left(-\nu_{2}/\sin(\theta)\right)}
%{\left|\sin(\theta)\right|}e^{-2\pi c \bm{i}(\sin(\theta)\nu_{1}+\cos(\theta)\nu_{2})
%/\sin(\theta)}\\
%\nonumber
%&&e^{2\pi \bm{i}
%\bm{\nu}^T\bm{r}_{-\theta}\x}\;d^2\bm{\nu}\\
\overset{(1)}{=}\nonumber 2^j(\pm 1) (-\bm{i}) \int_{\nu_{1,\min}}^{\nu_{1,\max}}
\tilde{H}_{j,u}(\cos(\theta)\nu_1,\sin(\theta)\nu_1)
e^{-2\pi \bm{i}c \nu_{1}}\;e^{2\pi \bm{i}
\nu_{1}(\cos(\theta)x_1+\sin(\theta)x_2)}\\
&&\times{\mathrm{sgn}}\left(\nu_1\right)\;d\nu_{1}+\rho_{8}\ +O\left(\frac{1}{N}\right)=  (\pm) U_{\bm{\xi}}^{(e)}+\rho_{8}+O\left(\frac{1}{N}\right),
\end{eqnarray}
where the value of $\pm$ depends on the value of $\theta.$
Similarly it transpires that
$
\mu_3^{(h,u)}=  \pm U_{\bm{\xi}}^{(e)}+\rho_{9}+O\left(\frac{1}{N}\right)$ and
$\mu_4^{(h,u)}= \pm W_{\bm{\xi}}^{(e)}+\rho_{10}+O\left(\frac{1}{N}\right).$
$\rho_8,$ $\rho_9$ and $\rho_{10}$ are constants depending on the variability
of $a_e(\x).$

\section*{C: Statistics of the Normal Vector}
For simplicity set $\Delta x =1$ when deriving the statistical properties
of the coefficients. $\epsilon_{\x}$ has a spectral representation:
$
\epsilon_{\x}=\int_{\bm{\Omega}} dZ_{\epsilon}(\f)e^{2\bm{i}\pi \f^T\x},$ where
$Z_{\epsilon}(\f)$ is a complex-valued orthogonal increment process, see
\cite{gikhman}[p.~244], {\em i.e.} $\E\left(dZ_{\epsilon}(\f)dZ_{\epsilon}^*(\f^{\prime})\right)=0$
if $\f \neq \f^{\prime}.$
The DWT is represented by subsampling
the MODWT:\\
$
\widetilde{{\mathbf{n}}}^{(s,u)}=
\left[\widetilde{W}_{\bm{\xi}_{\x}}^{(\epsilon,s,0)},\dots,
\widetilde{W}_{\bm{\xi}_{\x}}^{(\epsilon,s,L)}\right]^T,
$
noting that:
\begin{equation}
\label{covy1}
{\mathrm{cov}}\left(n^{(s,l_1,u)},
n^{(s,l_2,u)}\right)=2^{2j} {\mathrm{cov}}\left(
\tilde{n}^{(s,l_1,u)},
\tilde{n}^{(s,l_2,u)}\right).\end{equation}
%\subsection*{Stochastic Properties of $n_1^{(s,u)}$}
We have that
$
\tilde{n}^{(s,u)}_1=\int_{\bm{\Omega}} \tilde{H}_{j,u}\left(\f\right)\;dZ_{\epsilon}(\f)e^{2\pi i\f^T\x},$
$\tilde{H}_{j,1}\left(\f\right)=
\tilde{H}_j(f_1)\tilde{H}_j(f_2),$ $\tilde{H}_{j,2}\left(\f\right)=
\tilde{H}_j(f_1)\tilde{G}_j(f_2),$ $\tilde{H}_{j,3}\left(\f\right)=
\tilde{G}_j(f_1)\tilde{H}_j(f_2)$ and $\tilde{H}_{j,4}\left(\f\right)=
\tilde{G}_j(f_1)\tilde{G}_j(f_2).$
We approximate the magnitude of the wavelet filters as exact bandpass filters
- see for example Nielsen \cite{nielsen} for a discussion of such approximations,
and optimal filters to use. That is:
\begin{equation}
\left|\tilde{H}_{j}\left(f\right)\right|^2=\left\{
\begin{array}{lcr}
1 &{\mathrm{if}}& \left|f\right|\in \left[\frac{1}{2^{j+1}}
,\frac{1}{2^{j}}\right)\\
0 & {\mathrm{if}}& \left|f\right|\in \left[\frac{1}{2^{j+1}}
,\frac{1}{2^{j}}\right),
\end{array}\right.,\quad
\left|\tilde{G}_{j}\left(f\right)\right|^2=\left\{
\begin{array}{lcr}
1 &{\mathrm{if}}& f\in \left(-\frac{1}{2^{j+1}}
,\frac{1}{2^{j+1}}\right)\\
0 & {\mathrm{if}}& \left|f\right|\in \left[\frac{1}{2^{j+1}}
,\frac{1}{2}\right).
\end{array}\right. 
\label{GHapprox}
\end{equation}
%Thus we find that for $s=r$ and $s=h:$
\begin{eqnarray*}
{\mathrm{Var}}\left(n_{1}^{(s,u)}\right)%&=&2^{2j}
%E\left(
%\tilde{n}_{1}^{(s,u)}
%\tilde{n}_{1}^{(s,u)*}
%\right)\\ \nonumber
&=&2^{2j} E\left(\int\int_{\bm{\Omega}}  \int\int_{\bm{\Omega}}
\tilde{H}_{j,u}\left(\f\right) dZ_{\epsilon}(\f)e^{2\pi i (\f-\f^{\prime})^T\x }\tilde{H}_{j,u}^*\left(\f^{\prime}\right)
dZ_{\epsilon}^*(\f^{\prime})
\right)= \sigma^2%+O\left(\frac{1}{N}\right).
%&=&2^{2j}\sigma_1^2 \sigma_2^2\int_{-\frac{1}{2}}^{\frac{1}{2}} 
%\int_{-\frac{1}{2}}^{\frac{1}{2}}  \left|\tilde{H}_{j}\left(f_1\right)\right|^2
%\left|\tilde{H}_{j}\left(f_2\right)\right|^2\;df_1\;df_2\\
%&=& 4\sigma_1^2 \sigma_2^2\int_{\frac{1}{2^{j+1}}}^
%{\frac{1}{2^{j}}}\int_{\frac{1}{2^{j+1}}}^
%{\frac{1}{2^{j}}}\;df_1\;df_2\\
%&=&4\sigma_1^2 \sigma_2^2\left[s_1\right]_{\frac{1}{2^{j+1}}}^
%{\frac{1}{2^{j}}}\left[s_2\right]_{\frac{1}{2^{j+1}}}^
%{\frac{1}{2^{j}}}\\
%&=&4\sigma_1^2 \sigma_2^2\left[\frac{1}{2^{j}}-\frac{1}{2^{j+1}}\right]
%\left[\frac{1}{2^{j}}-\frac{1}{2^{j+1}}\right]\\
%&=&2^{2j+2}\sigma_1^2 \sigma_2^2\left(\frac{1}{2^{j+1}}\right)
%\left(\frac{1}{2^{j+1}}\right)\\
\end{eqnarray*}
%This then gives the variance of the first component corresponding to the
%wavelet transform of the observed image, in both the instance of the RTs
%and the HCTs signal. 

% varry2grejs3.m checks these calcs - they look aok
%From equation (\ref{covy1}) it follows that
\subsection*{Proof of Propositions \ref{proppy1} \& \ref{proppy2}}
\begin{eqnarray}
\nonumber
\var\left(n^{(s,u)}_l
\right)
%&=&
%2^{2j}\var\left(\tilde{n}^{(s,u)}_l
%\right)\\ \nonumber
&=& 2^{2j}
E\left(\int\int_{\bm{\Omega}}\int \int_{\bm{\Omega}}\tilde{H}_{j,u}(\f)\tilde{H}_{j,u}^*(\f^{\prime})
e^{2\pi \bm{i} (\f^T - \f^{\prime T})\x } \nonumber
dZ_{\breve{\epsilon}}(\f)
dZ_{\breve{\epsilon}}^*(\f^{\prime}) \right)\\ \nonumber
&=&2^{2j} \int_{\bm{\Omega}}
\left|V^{(s,l)}_D(\f)\right|^2\left|\tilde{H}_{j,u}(\f)\right|^2\;d^2\f \label{var}\\
&=&2^{2j}\int_{\bm{\Omega}}
\left|V^{(s,l)}(\f)\right|^2\left|\tilde{H}_{j,u}(\f)\right|^2\;d^2\f+O\left(
\frac{1}{N}\right)\equiv a^{(s,u)}_l+O\left(
\frac{1}{N}\right),
\nonumber
\end{eqnarray}
The latter defining $a^{(s,u)}_l,$ to be explicitly determined
for $s=r,\;h$ and $u=1,2,3,4.$ For $l_1\neq l_2,$ where $l_1,\;l_2\neq 0,$ we determine that:
\begin{eqnarray}
\nonumber
\cov\left(n^{(s,u)}_{l_1},
n^{(s,u)}_{l_2}
\right)
%&=& 2^{2j}\cov\left(\tilde{n}^{(s,u)}_{l_1},
%\tilde{n}^{(s,u)}_{l_2}
%\right)
%\\
%&=&E\left(\sum_{\bm{k}^{\prime}} \sum_{\bm{k}^{\prime\prime}}
%h_{\bm{\xi}^{\prime}}
%h_{\bm{\xi}^{\prime \prime}}^*
%\breve{\epsilon}_{2^j(\bm{k}^{\prime}+\bm{1})}^{(s,l_1)}
%\breve{\epsilon}_{2^j(\bm{k}^{\prime}+\bm{1})}^{(s,l_2)*}\right)\nonumber
%\\ \nonumber
&=& 2^{2j} \int_{\bm{\Omega}} \sigma^2 V^{(s,l_1)}_D(\f)
V_D^{(s,l_2)*}(\f)
\left|\tilde{H}_{j,u}(\f)\right|^2\;d^2\f \nonumber\\
&=&2^{2j}\sum\sum_{\x}\breve{h}_{\bm{\xi}_{\x}}^{(s,l_1)}\breve{h}_{\bm{\xi}_{\x}}^{(s,l_2)*}
+O\left(\frac{1}{N}\right)=O\left(\frac{1}{N}\right),
\label{covariance}
\end{eqnarray}
by property \ref{cond1} as $\tilde{h}_{\bm{\xi}_{\x}}$ is separable.
Therefore 
\begin{eqnarray}
\nonumber
\E\left(\sum_{l=1}^{L} n_{l}^{(s,u)2}
\right)
&=&2^{2j}\int_{\bm{\Omega}}\sum_{l=1}^L \left|V^{(s,l)}(\f) \right|^2 \sigma^2 \left|H_{j,u}(\f)\right|^2\;d^2\f+O\left(\frac{1}{N}\right)\\
&=&2^{2j}\int_{\bm{\Omega}}C_{L}^{(s)} \nonumber
\sigma^2 \left|H_{j,u}(\f)\right|^2\;d^2\f+O\left(\frac{1}{N}\right)\\
&=&C_L^{(s)} \var\left(n_{1}^{(s,u)2} \right)+O\left(\frac{1}{N}\right).
\end{eqnarray}
Hence the total energy of the noise associated with the total magnitude square
of the added quadrature components is a constant times the variance of the
original signal.
%
%\subsection*{Proof of Proposition \ref{proppy2}}
%Also we determine that:
\begin{eqnarray}
\cov\left(n_{l}^{(s,u)} n_1^{(s,u)}
\right)\nonumber
&=&2^{2j}E\left(\int\int_{\bm{\Omega}}
\int \int_{\bm{\Omega}}H_{j,u}(\f)H_{j,u}^*(\f^{\prime})
dZ_{\breve{\epsilon}}(\f)\right.\\
&& \nonumber \left.
dZ_{\epsilon}^*(\f^{\prime}) e^{2\pi \bm{i}(\f-\f^{\prime})\x} \right)
=\int_{\bm{\Omega}} \sigma^2 V^{(s,l)}(\f) H_{j,u}(\f)
H_{j,u}^*(\f)\;d^2\f\\
&&+O\left(\frac{1}{N}\right)\nonumber =\sum_{\x} \breve{h}_{\bm{\xi}_{\x}}^{(s,l)}
h_{\bm{\xi}_{\x}}^{(s,l)}+O\left(\frac{1}{N}\right)
=O\left(\frac{1}{N}\right).
\end{eqnarray}

\subsection*{Proof of Lemma \ref{rieszcoef}}
%Given it was assumed that the noise was Gaussian, to determine the relevant
%distributions we
%wish to determine the joint second order structure of the observed
%image and the additional components observed, and especially determine the
%variance of each component. We find:
Given the noise was Gaussian and zero-mean we only need to determine the
second order structure to deduce the Lemma, using equation (\ref{covy1}).
\begin{eqnarray*}
{\mathrm{Var}}\left(n_2^{(r,1)}\right)%&=&2^{2j}
%E\left(
%\tilde{n}_2^{(r,1)}
%\tilde{n}_2^{(r,1)*}
%\right)\\
&=&2^{2j} E\left(\int\int_{\bm{\Omega}}
\int\int_{\bm{\Omega}}
 \tilde{H}_{j}\left(f_1\right)
\tilde{H}_{j}\left(f_2\right)\;dZ_{\epsilon}(\f) e^{2\pi \bm{i}(\f-\f^{\prime})^T\x}\right.\\
&&\left.\frac{f_1}{\sqrt{f_1^2+f_2^2}}\tilde{H}_{j}^*\left(f_1^{\prime}\right)
\tilde{H}_{j}^*\left(f_2^{\prime}
\right)\;dZ_{\epsilon}(\f^{\prime})\frac{f_1^{\prime}}
{\sqrt{f_1^{\prime 2}+f_2^{\prime 2}}}\right)+O\left(\frac{1}{N}\right)\\
&=&2^{2j}\sigma^2\int\int_{\bm{\Omega}}  \frac{f_1^2}{
f_1^2+f_2^2}\left|\tilde{H}_{j}\left(f_1\right)\right|^2
\left|\tilde{H}_{j}\left(f_2\right)\right|^2\;d^2\f+O\left(\frac{1}{N}\right).
\end{eqnarray*}
%Clearly by symmetry it follows that
\begin{eqnarray*}
{\mathrm{Var}}\left(n_3^{(r,1)}\right)%&=&%2^{2j}
%E\left(
%\tilde{n}_2^{(r,1)}\tilde{n}_2^{(r,1)*}
%\right)\\
&=&2^{2j}\sigma^2\int\int_{\bm{\Omega}}  \frac{f_2^2}{
f_1^2+f_2^2}\left|\tilde{H}_{j}\left(f_1\right)\right|^2
\left|\tilde{H}_{j}\left(f_2\right)\right|^2\;d^2\f+O\left(\frac{1}{N}\right).
\end{eqnarray*}
Hence it follows that
\begin{equation}
{\mathrm{Var}}\left(n_2^{(r,1)}\right)+
{\mathrm{Var}}\left(n_3^{(r,1)}\right)=
{\mathrm{Var}}\left(n_1^{(r,1)}\right),
\end{equation}
and the two variances are obviously equal. Also:
\begin{equation}
{\mathrm{Var}}\left(n_2^{(r,u)}\right)+
{\mathrm{Var}}\left(n_3^{(r,u)}\right)=
{\mathrm{Var}}\left(n_1^{(r,u)}\right),\;u=1,2,3,4,
\end{equation}
\begin{eqnarray*}
{\mathrm{Var}}\left(n_2^{(r,1)}\right)
&=&2^{2j+2}\sigma^2\int_{\frac{1}{2^{j+1}}}^
{\frac{1}{2^{j}}} \int_{\frac{1}{2^{j+1}}}^
{\frac{1}{2^{j}}}  \frac{f_1^2}{
f_1^2+f_2^2}\;df_2\;df_1+O\left(\frac{1}{N}\right)\\
&=&2^{2j+2}\sigma^2\int_{\frac{1}{2^{j+1}}}^
{\frac{1}{2^{j}}}f_1 \left[\tan^{-1}\left(\frac{f_2}{f_1}\right)\right]_{f_2=\frac{1}{2^{j+1}}}
^{\frac{1}{2^{j}}} \;df_1+O\left(\frac{1}{N} \right)\\
%&=&4\sigma_1^2 \sigma_2^2\int_{\frac{1}{2^{j+1}}}^
%{\frac{1}{2^{j}}}f_1 \left(\tan^{-1}\left(\frac{1}{2^{j}f_1}\right)
%-\tan^{-1}\left(\frac{1}{2^{j+1}f_1}\right)\right) \;df_1\\
%&=&4\sigma_1^2 \sigma_2^2 2^{-2j}\int_{\frac{1}{2}}^
%{1}f_1 \tan^{-1}\left(\frac{1}{f_1}\right)\;df_1\\
%&&-4\sigma_1^2 \sigma_2^2 2^{-2j-2}\int_{1}^
%{2}f_1 \tan^{-1}\left(\frac{1}{f_1}\right)\;df_1\\
%&=&\sigma_1^2 \sigma_2^2 2^{-2j+2}
%\left[
%\frac{1}{2}f_1^2\tan^{-1}\left(\frac{1}{f_1}\right)+\frac{1}{2}f_1
%-\frac{1}{2}\tan^{-1}\left(f_1\right)
%\right]_{\frac{1}{2}}^{1}\\
%&&-\sigma_1^2 \sigma_2^2 2^{-2j}
%\left[
%\frac{1}{2}f_1^2\tan^{-1}\left(\frac{1}{f_1}\right)+\frac{1}{2}f_1
%-\frac{1}{2}\tan^{-1}\left(f_1\right)
%\right]_{1}^{2}
%\\
%&=&\sigma_1^2 \sigma_2^2 2^{-2j+2}
%\left[\frac{1}{2}-
%\frac{1}{8}\tan^{-1}\left(2\right)-\frac{1}{4}
%+\frac{1}{2}\tan^{-1}\left(\frac{1}{2}\right)
%\right]\\
%&&-\sigma_1^2 \sigma_2^2 2^{-2j}
%\left[
%\frac{2^2}{2}\tan^{-1}\left(\frac{1}{2}\right)+1
%-\frac{1}{2}\tan^{-1}\left(2\right)-\frac{1}{2}
%\right]
%\\
&=&\sigma^2 \left(\frac{1}{4}
-\frac{1}{8}\tan^{-1}\left(2\right)+\frac{1}{2}\tan^{-1}\left(\frac{1}{2}\right)-
\left(\frac{1}{8}
-\frac{1}{8}\tan^{-1}\left(2\right)+\frac{1}{2}\tan^{-1}\left(\frac{1}{2}\right)
\right)
\right)\\
&& \nonumber+O\left(\frac{1}{N} \right)= \frac{\sigma^2}{2}+O\left(\frac{1}{N} \right)={\mathrm{Var}}\left(n_3^{(r,u)}\right)=\sigma^2 a^{(r,1)},
\end{eqnarray*}
\begin{eqnarray*}
{\mathrm{Var}}\left(n_2^{(r,2)}\right)
&=&2^{2j+2}\sigma^2\int_{\frac{1}{2^{j+1}}}^
{\frac{1}{2^{j}}} \int_{0}^
{\frac{1}{2^{j+1}}}  \frac{f_1^2}{
f_1^2+f_2^2}\;df_2\;df_1+O\left(\frac{1}{N} \right)
%&=&4\sigma_1^2 \sigma_2^2\int_{\frac{1}{2^{j+1}}}^
%{\frac{1}{2^{j}}} f_1\left[\tan^{-1}\left(\frac{f_2}{f_1}\right)
%\right]^{\frac{1}{2^{j+1}}}_{f_2=0}\;df_1\\
%&=&4\sigma_1^2 \sigma_2^2\int_{\frac{1}{2^{j+1}}}^
%{\frac{1}{2^{j}}} f_1\tan^{-1}\left(\frac{1}{2^{j+1} f_1}\right)
%\;df_1\\
%&=&4\sigma_1^2 \sigma_2^2 2^{-2(j+1)}\int_{1}^
%{2} f_1\tan^{-1}\left(\frac{1}{f_1}\right)
%\;df_1\\
%&=&\sigma_1^2 \sigma_2^2 2^{-2j}\left[
%\frac{1}{2}f_1^2\tan^{-1}\left(\frac{1}{f_1}\right)+\frac{1}{2}f_1
%-\frac{1}{2}\tan^{-1}\left(f_1\right)
%\right]_{1}^{2}\\
%&=&\sigma_1^2 \sigma_2^2 2^{-2j}\left(
%2\tan^{-1}\left(\frac{1}{2}\right)+1
%-\frac{1}{2}\tan^{-1}\left(2\right)-\frac{1}{2}
%\right)\\
=\sigma^2 \left(\frac{1}{2}+
2\tan^{-1}\left(\frac{1}{2}\right)\right.\\
&&\left.
-\frac{1}{2}\tan^{-1}\left(2\right)
\right)
+O\left(\frac{1}{N}\right)=\sigma^2 a^{(r,2)}+O\left(\frac{1}{N}\right)\approx 0.8737\sigma^2 .
\end{eqnarray*}
Finally note that:
\begin{eqnarray*}
{\mathrm{Var}}\left(n_2^{(r,3)}\right)
&=&2^{2j+2}\sigma^2
\int_{0}^{\frac{1}{2^{j+1}}}
\int_{\frac{1}{2^{j+1}}}^
{\frac{1}{2^{j}}}   \frac{f_1^2}{
f_1^2+f_2^2}\;df_2\;df_1+O\left(\frac{1}{N}\right)=
%&=&4\sigma_1^2 \sigma_2^2
%\int_{0}^{\frac{1}{2^{j+1}}}
%\int_{\frac{1}{2^{j+1}}}^
%{\frac{1}{2^{j}}}   \left(1-\frac{f_2^2}{
%f_1^2+f_2^2}\right)\;df_2\;df_1\\
%&=&\sigma_1^2 \sigma_2^2 2^{-2j}-\sigma_1^2 \sigma_2^2 2^{-2j}a_j^{(2)}\\
\sigma^2 \left(1-a^{(r,2)}\right)+O\left(\frac{1}{N}\right)\\
&\equiv &\sigma^2  a^{(r,3)}\approx 0.1263\sigma^2,
\end{eqnarray*}
\begin{eqnarray*}
{\mathrm{Var}}\left(n_{j,k_1,k_2}^{(r,1,4)}\right)
&=&2^{2j+2}\sigma^2
\int_{0}^
{\frac{1}{2^{j+1}}} 
\int_{0}^
{\frac{1}{2^{j+1}}}  \frac{f_1^2}{
f_1^2+f_2^2}\;df_2\;df_1
%&=&4\sigma_1^2 \sigma_2^2\int_{0}^
%{\frac{1}{2^{j+1}}} f_1\left[\tan^{-1}\left(\frac{f_2}{f_1}\right)
%\right]^{\frac{1}{2^{j+1}}}_{f_2=0}\;df_1\\
%&=&4\sigma_1^2 \sigma_2^2\int_{0}^
%{\frac{1}{2^{j+1}}} f_1\tan^{-1}\left(\frac{1}{2^{j+1} f_1}\right)
%\;df_1\\
%&=&4\sigma_1^2 \sigma_2^2 2^{-2(j+1)} \int_{0}^
%{1} f_1\tan^{-1}\left(\frac{1}{ f_1}\right)
%\;df_1\\
%&=&\sigma_1^2 \sigma_2^2 2^{-2j}
%\left[
%\frac{1}{2}f_1^2\tan^{-1}\left(\frac{1}{f_1}\right)+\frac{1}{2}f_1
%-\frac{1}{2}\tan^{-1}\left(f_1\right)
%\right]_{0}^{1}
%\\
=\frac{\sigma^2}{2}+O\left(\frac{1}{N}\right)=
\sigma^2 a^{(r,4)}+O\left(\frac{1}{N}\right).
\end{eqnarray*}
Clearly we can find the variance of the second RT by permuting
the order of the spatial variable in the integration, and this then completes the variance calculations.
From the proofs of propositions (\ref{proppy1}) and (\ref{proppy2}) we can note that the components of ${\mathbf{n}}^{(r,u)}$
for $u=1,\;2,\;3,\;4$ are uncorrelated up to $O\left(\frac{1}{N}\right),$
and this can also be shown by direct calculation, {\em mutatis mutandis} the
calculations given above. Thus as $\epsilon_{\x}$ was zero-mean Gaussian, and we are
forming linear combinations to obtain ${\mathbf{n}}^{(r,u)},$ the stated
result follows from the expressions for the covariances of the components.

\subsection*{Proof of Lemma \ref{rieszmag}}
a) First consider $u=1,4$ so that the variance of the two Riesz components is
$1/2.$ Then by Lemma \ref{rieszcoef} it follows directly that
$(C_L^{(s)}+1)M_{\bm{\xi}}^{(\epsilon,r)2} 
\overset{\cal L}{=}\bm{Z}^{(r,u)T}\bm{Z}^{(r,u)}+O\left(
\frac{1}{N}\right),\;T_1=\bm{Z}^{(r,u)T}\bm{Z}^{(r,u)}\sim
\chi^2_1+\frac{1}{2} \chi^2_1+\frac{1}{2}\chi^2_1.$
$T_1$ has a Moment Generating Function (MGF) given by
$M_{T_1}(s)=\frac{1}{\sqrt{1-2s}}\frac{1}{1-s},$
%and we shall consider the inversion to the pdf. Give the Laplace transform
%of a density the notation ${\mathcal{L}}\left\{\cdot\right\}.$
and thus $f_{T_1}(t)= e^{-t}\frac{2}{\sqrt{\pi}}\int_{0}^{\sqrt{\frac{t}{2}}}
e^{u^2}\;du.$
%\begin{eqnarray*}
%{\mathcal{L}}^{-1}\left\{
%M_{T}(-s)\right\}\left(t\right)
%&=&{\mathcal{L}}^{-1}\left\{
%\frac{1}{\sqrt{1+2s}}\frac{1}{1+s}\right\}\left(t\right)\\
%&=&\frac{1}{2i\pi}\int_{\gamma-i\infty}^{\gamma+i\infty}
%M_{T}(-x)e^{t s}\;ds
%\\
%&=&\sqrt{\frac{2t}{\pi}}e^{-\frac{t}{2}}\;_1 F_1\left(1,\frac{3}{2},-\frac{t}{2}\right)
%\\
%&=&\sqrt{\frac{2t}{\pi}}e^{-\frac{t}{2}}\frac{1}{2}\sqrt{\frac{2\pi}{t}}
%e^{-\frac{t}{2}}{\mathrm{erfi}}\left(\sqrt{\frac{t}{2}}\right)
%\\
%&=&e^{-t}{\mathrm{erfi}}\left(\sqrt{\frac{t}{2}}\right)\\
%&=& e^{-t}\frac{2}{\sqrt{\pi}}\int_{0}^{\sqrt{\frac{t}{2}}} 
%e^{u^2}\;du\\
%&=&
%%f_{T}(t).
%\end{eqnarray*}
The probability that $T_1$ does not exceed
$\lambda^2$ is given by 
\begin{eqnarray}
\nonumber
P\left(T_1<\lambda^2
\right)&=&%\int_0^{\lambda^2}f_{T}(x)\;dx\\
%\nonumber
%&=&\int_0^{\lambda^2}\int_{0}^{\sqrt{\frac{x}{2}}}e^{-x}\frac{2}{\sqrt{\pi}}
%e^{u^2}\;du\;dx\\
%\nonumber
%&=&1-\int_{\lambda^2}^{\infty}\int_{0}^{\sqrt{\frac{x}{2}}}e^{-x}\frac{2}{\sqrt{\pi}}
%e^{u^2}\;du\;dx\\
%\nonumber
%&=&1+\left[e^{-x}
%\int_{0}^{\sqrt{\frac{x}{2}}}
%\frac{2}{\sqrt{\pi}}
%e^{u^2}\;du\right]_{\lambda^2}^{\infty}
%-\int_{\lambda^2}^{\infty}e^{-x}
%\frac{2}{\sqrt{\pi}}
%e^{\frac{x}{2}}\frac{1}{\sqrt{2}}\frac{1}{2\sqrt{x}}\;dx\\
%\nonumber
%&=&1+\left[e^{-x}
%\int_{0}^{\sqrt{\frac{x}{2}}}
%\frac{2}{\sqrt{\pi}}
%e^{u^2}\;du\right]_{\lambda^2}^{\infty}
%-\int_{\lambda^2}^{\infty}
%\frac{1}{\sqrt{2\pi }} x^{-\frac{1}{2}}
%e^{-\frac{x}{2}}\;dx\\
%\nonumber
%&=&
1+\left[e^{-x}
\int_{0}^{\sqrt{\frac{x}{2}}}
\frac{2}{\sqrt{\pi}}
e^{u^2}\;du\right]_{\lambda^2}^{\infty}
-\int_{\lambda^2}^{\infty}
\frac{1}{\sqrt{2\pi }} x^{-\frac{1}{2}}
e^{-\frac{x}{2}}\;dx\\
&=&1-\sqrt{\frac{8}{\pi}}\frac{1}{\lambda}e^{-\lambda^2/2}+
e^{-\lambda^2/2}\left(O\left(\frac{1}{\lambda^3}\right)\right)
\label{asprob},
\end{eqnarray}
where the cdf be found in Grad \& Solomon \cite{grad}[p.~472], and the function
is expanded as $\lambda\rightarrow \infty.$

b) Consider now $u=2,\;3.$
Wlog assume that $a^{(r,u)}>1/2,$ and otherwise relabel $a^{(r,u)}$ and $1-a^{(r,u)}$
suitably. If $a^{(r,u)}=1-a^{(r,u)}$ this collapses to the case given above.
Define $T^{\prime}=T_1/2$ so that $T^{\prime}=\sum_{i=1}^3 a_i X_i^2$
where the $X_i$ are iid Gaussian random variates with zero mean and unit
variance, where $\sum_{i=1}^3 a_i=1.$ Then the MGF
is by Grad \& Solomon \cite{grad}
\begin{eqnarray*}
M_{T^{\prime}}(s)=
M_{T_1}(s/2)
=\left(1-2(1-a^{(r,u)})s/2\right)^{-\frac{1}{2}}
\left(1-2a^{(r,u)}s/2\right)^{-\frac{1}{2}}
\left(1-2s/2\right)^{-\frac{1}{2}}
=\prod_{l=1}^3\left(1-2a_l s\right)^{-\frac{1}{2}},
\end{eqnarray*}
with $a_1=(1-a^{(r,u)})/2,\;a_2=a^{(r,u)}/2,\;a_3=1/2.$
Defining
$c_1=2/(1-a^{(r,u)}),\;c_2=2/a^{(r,u)},\;c_3=2,$
and hence as we assumed $1>a^{(r,u)}>1-a^{(r,u)}$ we may note that $c_1\ge c_2 \ge c_3$ and thus
agrees with \cite{grad}'s notation. For future reference note that
$c_1+c_2=2/(a^{(r,u)}(1-a^{(r,u)})),\;u=1,\dots,4,$ and
we may rewrite 
$M_{T^{\prime}}(s)=\prod_{j=1}^{3} \left(1-2s/c_j\right)^{-1/2}.$
Using results from Grad and Solomon we may determine:
\begin{eqnarray}
\nonumber
F_{T_1}\left(t\right)=
F_{T^{\prime}}\left(t/2\right)
\nonumber
%&=&\left(1-h_1(t/2)-h_2(t/2)\right)\\
%\nonumber
%&=&\left(1
%-\frac{\sqrt{c_2}}{\sqrt{\pi}}
%e^{-c_3t/4}\frac{2}{\sqrt{t}}
%\left[\frac{\sqrt{ 2c_3}e^{-(c_1-c_3)t/4}}{\sqrt{c_1(c_2-c_1)(c_3-c_1)}}
%+\frac{\sqrt{c_1 }}{\sqrt{c_3(c_3-c_1)(c_3-c_2)}}\right]\right)
%+o(1)
%\nonumber\\
=\left(1-A e^{-c_3t/4}\frac{1}{\sqrt{t}}\right)+o(1).
\end{eqnarray}
For suitably defined constant $A.$ From these formulae we can thus consider the probability of an observation exceeding a large threshold,
which will be necessary for the selection of an appropriate threshold.

\subsection*{Proof of Lemma \ref{hypercoef}}
We establish 
\begin{eqnarray*}
{\mathrm{Var}}\left(n_{2}^{(h,1)}\right)
%&=&
%2^{2j}E\left(
%\tilde{n}_{2}^{(h,1)}
%\tilde{n}_{2}^{(h,1)*}
%\right)
%&=&2^{2j}E\left(\int_{-\frac{1}{2}}^{\frac{1}{2}} 
%\int_{-\frac{1}{2}}^{\frac{1}{2}} \int_{-\frac{1}{2}}^{\frac{1}{2}} 
%\int_{-\frac{1}{2}}^{\frac{1}{2}} \tilde{H}_{j}\left(f_1\right)
%\tilde{H}_{j}\left(f_2\right)\;dZ(f_1)\;dZ(f_2)\right. \\
%&&\left. e^{2\pi i(f_1 x_1+
%f_2 x_2)}{\mathrm{sgn}}\left(f_1\right)\tilde{H}_{j}^*\left(f_3\right)
%\tilde{H}_{j}^*\left(f_4\right)\right.\\
%&&\left.\;dZ^*(f_3)\;dZ^*(f_4)e^{-2\pi i(f_3 x_1+
%f_4 x_2)}{\mathrm{sgn}}\left(f_3\right)\right)\\
%&=&2^{2j}\sigma_1^2 \sigma_2^2\int_{-\frac{1}{2}}^{\frac{1}{2}} 
%\int_{-\frac{1}{2}}^{\frac{1}{2}} \left|\tilde{H}_{j}\left(f_1\right)\right|^2
%\left|\tilde{H}_{j}\left(f_2\right)\right|^2\;df_1\;df_2\\
=\sigma^2+O\left(\frac{1}{N}\right)={\mathrm{Var}}\left(n_{l}^{(h,u)}\right),\;u=1,2,3,4,\;l=1,2,3,4.
\end{eqnarray*}
These results follow trivially from the form of the partial HT \cite{hahn1996}.
From the proofs of propositions (\ref{proppy1}) and (\ref{proppy2}) we can note that the components of ${\mathbf{n}}^{(h,u)}$
for a fixed value of $u=1,\;2,\;3,\;4$ are uncorrelated up to $O\left(\frac{1}{N}\right),$
and this can also be shown by direct calculation {\em mutatis mutandis} the
calculations given above.

\subsection*{Proof of Lemma \ref{rieszmag}}
$K$ wavelet coefficients will be thresholded where
$K_1$ magnitudes have the distribution given when $u=1,4$ and $K_2$ have
the distribution that follows from $u=2,3,$ where $K_1+K_2=K,$ and 
$K_1,\;K_2,\;K_3=O(K),$ so that by Dykstra \cite{dykstra}
\begin{eqnarray}
\nonumber
P\left({\mathcal{M}}_K< \lambda_K^2\right)
&=& \left(1-\sqrt{\frac{8}{\pi}}\frac{1}{\lambda_K}e^{-\lambda_K^2/2}\right)^{K_1}
\left(1-A e^{-c_3\lambda_K^2/4}\frac{1}{\lambda_K}\right)^{K_2}
=\left[1-\frac{A_2}{\lambda_K}e^{-\lambda_K^2/2}
\right]^{K_3},
\label{thresh2}%\\
%\nonumber
%&=&\left[\left(1-\sqrt{\frac{8}{\pi}}\frac{1}{\sqrt{t}}e^{-t/2}\right)
%\left(1-A e^{-t/2}\frac{1}{\sqrt{t}}\right)\right]^{K/2}\\
%\nonumber
%&=&\left[1-\left(\sqrt{\frac{8}{\pi}} +A\right)\frac{1}{\sqrt{t}}e^{-t/2}
%+\sqrt{\frac{8}{\pi}}A\frac{1}{t}e^{-t}\right]^{K/2}\\
%\nonumber
%&=&\left[1-\frac{1}{\sqrt{t}}e^{-t/2}
%\left(\sqrt{\frac{8}{\pi}} +A-
%\sqrt{\frac{8}{\pi}}A\frac{1}{\sqrt{t}}e^{-t/2}\right)\right]^{K/2}\\
\end{eqnarray}
ignoring $o(1)$ terms in $K$
for suitably chosen constant $A_2,$ if $\lambda_K=O(\log[K]).$
Thus with
$
\lambda^2_K=2\log\left[K\right]+C_2
\log\left[\log\left[K\right]\right],
$
\begin{eqnarray*}
P\left({\mathcal{M}}_K< \lambda^2_K\right)&=&
\left(1-A_2 e^{-\lambda_K^2/2}\frac{1}{\lambda_K}\right)^{K_3}+o(1)\\
&=&
\left(1-A_2 e^{-(2\log\left[K\right]+C_2
\log\left[\log\left[K\right]\right])/2}\frac{1}{\sqrt{2\log\left[K\right]+C_2
\log\left[\log\left[K\right]\right]}}\right)^{K_3}
%&=&\left(1-A_2 \frac{1}{K\log(K)^{C_2/2}}
%\frac{1}{\sqrt{2\log\left[K\right]}}\right)^{K/2}
%\\
%&\rightarrow& e^{-A_2 \frac{2}{\log(K)^{C_2/2}}
%\frac{1}{\sqrt{2\log\left[K\right]}}}\\
\rightarrow  1,
\end{eqnarray*}
%\begin{equation}
%\label{require}
if $C_2/2>-1/2,\;\;i.e.\;\;C_2>-1,$ and
%\end{equation}
we take $C_2=0.$
\section*{D: Risk Calculations}
\subsection*{Proof of Theorem \ref{apy1}}
We firstly note from \cite{Marron}[p.~293] that the risk of regular hard
thresholding is given by
\begin{eqnarray}
\nonumber
R_{\theta}^{(c)}(\lambda)%&=&\int_{w^2\ge \lambda^2}
%(w-\theta_1)^2\phi(w-\theta_1)\;dw+\theta_1^2\int_{w^2<\lambda^2}\phi(w-\theta_1)\;dw\\
%\nonumber
&=&\int_{(w+\theta_1)^2\ge \lambda^2}
w^2\phi(w)\;dw+\theta_1^2\int_{(w+\theta_1)^2 <\lambda^2}\phi(w)\;dw\\
&=&1+\int_{(w+\theta_1)^2 <\lambda^2}\left[\theta_1^2-w^2 \right]\phi(w)\;dw.
\end{eqnarray}
We may then note that the risk for the hyperanalytic threshold with
$\var\left(n_l^{(s,u)}\right)=\sigma^2 \sigma_l^2,$ for $l=0,\dots,L,$
where $\sigma_l^2$ takes the value $1$ or $a^{(s,u)},$
is given by:
%\begin{figure*}[t]
%\centerline{
%\includegraphics[scale=1]{risk1plot3.eps}
%\includegraphics[scale=1]{risk1plot4.eps}
%}
%\caption{\label{exx4}
%The risk associated with a thresholded coefficient using standard hard thresholding
%(solid line), the Riesz threshold (dotted line) or the Hypercomplex threshold
%(dash dotted line).
%}
%\end{figure*}
\begin{eqnarray}
\nonumber
R_{\theta}^{(s)}(\lambda)&=&
\theta_1^2\int_{\sum_l w_l^2 \le  \lambda^2}\prod_{l}\sigma^{-1}_l
\phi\left(\frac{w_l-\theta_l}{\sigma_l}\right)\;dw_l+
\int_{\sum_l w_l^2 > \lambda^2}\left[w_1-\theta_1 \right]^2\prod_{l}\sigma^{-1}_l
\phi\left(\frac{w_l-\theta_l}{\sigma_l}\right)\;dw_l\\
&=&1+\int_{\sum_l (w_l+\theta_l)^2 \le  \lambda^2}\left[\theta_1^2-w_1^2 \right]\prod_{l}\sigma^{-1}_l
\phi\left(\frac{w_l}{\sigma_l}\right)\;dw_l.
\end{eqnarray}
\subsection*{Proof of Corollary \ref{risk0}}
%
%In this case we get
\begin{eqnarray}
%\nonumber
R_{0}^{(c)}(\lambda)&=&
1-\int_{-\lambda}^{\lambda} w^2\frac{1}{\sqrt{2\pi}}e^{-\frac{1}{2}w^2}\;dw
%=1-\int_0^{\frac{1}{2}\lambda^2} \frac{u^{1/2}}{\sqrt{2\pi}}e^{-u}\;du
\label{riskhard}
=1-\gamma\left(1/2,\frac{1}{2}\lambda^2\right)=e^{-\frac{1}{2}\lambda^2}
\left(\frac{\sqrt{2}}{\sqrt{\pi}\lambda}+O\left(\lambda^{-3}\right)\right).
\end{eqnarray}
%The risk of the `analytic' threshold is given by:
\begin{eqnarray}
\nonumber
R_{0}^{(a)}(\lambda)&=&
1-\int_{w_1^2+w_2^2 \le \lambda^2} w_1^2\frac{1}{2\pi}e^{-\frac{1}{2}(w_1+w_2)^2}\;dw_1\;dw_2
%=1-\frac{1}{2\pi}\int_{0}^{\lambda}\int_0^{2\pi} r^3 \cos^2(\theta) e^{-\frac{1}{2}r^2}\;dr\;d\theta\\
\label{hypertyp}
%&=&1-\int_{0}^{\frac{1}{2}\lambda^2}u  e^{-u}\;du
%=1-\gamma(2,\frac{1}{2}\lambda^2)
=e^{-\frac{1}{2}\lambda^2}\left(1+\frac{1}{2}\lambda^2\right).
\end{eqnarray}
Furthermore the risk at $\bm{\theta}=\bm{0}$ can also be found for the other
hyperanalytic thresholds. We note that for $s=r$ and $u=1,\;4,$ denoted by
$s=r_1$ in the Figures:
\begin{eqnarray}
\nonumber
R_{0}^{(r_1)}(\lambda)&=&
1-\int_{w_1^2+w_2^2+w_3^2\le \lambda^2} w_1^2\frac{2}{(\sqrt{2\pi})^3}
e^{-\frac{1}{2}\left(w_1^2+2\left(
w_2^2+w_3^2\right)\right)}\;dw_1\; dw_2\; dw_3\\
\label{rieszriskatnought}
&=&1-\int_{-\lambda}^{\lambda}w_1^2\frac{2}{(\sqrt{2\pi})^3} e^{-\frac{1}{2}w_1^2}
\left[\int \int_{w_2^2+w_3^2\le \lambda^2-w_1^2} e^{-(w_2^2+w_3^2)}\; dw_2\; dw_3\right]\;dw_1\\
%\nonumber
%&=&1-\int_{-\lambda}^{\lambda}w_1^2\frac{2}{(\sqrt{2\pi})^3} e^{-\frac{1}{2}w_1^2}
%(\pi)\left[
%1-e^{-(\lambda^2-w_1^2)}
%\right]\;dw_1\\
%\nonumber
%&=&1-\int_{-\lambda}^{\lambda}w_1^2\frac{1}{(\sqrt{2\pi})} 
%\left[e^{-\frac{1}{2}w_1^2}-e^{-(\lambda^2-1/2 w_1^2)}
%\right]\;dw_1\\
%\nonumber
%&=&1-\gamma\left(1/2,\frac{1}{2}\lambda^2\right)+2
%\int_{\lambda}^{\sqrt{2} \lambda}(2\lambda^2-w_1^2)\frac{1}{(\sqrt{2\pi})}
%e^{-\frac{1}{2}w_1^2}\;dw_1\\
%\nonumber
%&=&1-\gamma\left(1/2,\frac{1}{2}\lambda^2\right)+4\lambda^2
%\left(\Phi(\sqrt{2} \lambda)-\Phi( \lambda)\right)-2\left(
%\gamma\left(1/2,\lambda^2\right)-
%\gamma\left(1/2,\frac{1}{2}\lambda^2\right)
%\right)\\ 
\nonumber
&=&1+\gamma\left(1/2,\frac{1}{2}\lambda^2\right)+4\lambda^2
\left(\Phi(\sqrt{2} \lambda)-\Phi( \lambda)\right)-2
\gamma\left(1/2,\lambda^2\right)=e^{-\frac{1}{2}\lambda^2}\left(\frac{4\sqrt{2}}{\sqrt{\pi}\lambda}+
O\left(\lambda^{-3} \right) \right),
\end{eqnarray}
whilst for $u=2,\;3$ denoted by $s=r_2,$ 
\begin{eqnarray}
\nonumber
R_{0}^{(r_2)}(\lambda)&=&
1-\int_{\sum_l w_l^2\le \lambda^2} w_1^2\frac{1}{(\sqrt{2\pi})^3\sqrt{a^{(r,u)}
(1-a^{(r,u)})}}
e^{-\frac{1}{2}\left(w_1^2+
w_2^2/a^{(r,u)}+w_3^2/(1-a^{(r,u)})\right)}\;dw_1\; dw_2\; dw_3%\\
%\nonumber
%&=&1-\int_{-\lambda}^{\lambda} w_1^2\frac{1}{(\sqrt{2\pi})^3\sqrt{a^{(u)}
%(1-a^{(u)})}}e^{-\frac{1}{2}w_1^2}\int_0^{\lambda^2-w_1^2} \int_0^{2\pi}
%e^{-\frac{r^2}{2}\left(
%\cos^2(\theta)/a^{(u)}+\sin^2(\theta)/(1-a^{(u)})\right)}\; 
%rdr\;d\theta  \;dw_1\\
%&=&1-\int_{-\lambda}^{\lambda} w_1^2\frac{1}{(\sqrt{2\pi})^3\sqrt{a^{(u)}
%(1-a^{(u)})}}e^{-\frac{1}{2}w_1^2}\int_0^{2\pi}\left[
%-\frac{e^{-\frac{r^2}{2}\left(
%\cos^2(\theta)/a^{(u)}+\sin^2(\theta)/(1-a^{(u)})\right)}}
%{\cos^2(\theta)/a^{(u)}+\sin^2(\theta)/(1-a^{(u)})}\right]\; d\theta  \;dw_1.
\end{eqnarray}
%and
\begin{eqnarray}
\nonumber
R_{0}^{(h)}(\lambda)&=&
1-\int_{w_1^2+w_2^2+w_3^2+w_4^2\le \lambda^2} w_1^2\frac{1}{(\sqrt{2\pi})^4}e^{-\frac{1}{2}\sum
w_l^2}\;dw\\
\label{hcriskatnought}
&=&1-\int_{0}^{\lambda} \int_{0}^{2\pi} \int_{0}^{\pi} \int_{0}^{\pi} 
r^2 \cos^2(\theta) 
\frac{1}{(\sqrt{2\pi})^4}e^{-\frac{1}{2} r^2} r^3 \sin^2(\theta)\sin (\phi)\;
dr\;d\theta\;d\phi\;d\lambda\\
\nonumber
&=& 1-\frac{(\frac{\pi}{4})(2\pi)}{4\pi^2}\int_{0}^{\lambda} 
r^5 
e^{-\frac{1}{2} r^2} 
dr =1-\frac{1}{8}\int_{0}^{\frac{1}{2}\lambda^2} 
4 s^2
e^{-s} 
ds=
e^{-\frac{1}{2}\lambda^2}\left(1+\frac{1}{2}\lambda^2+\frac{1}{8}\lambda^4\right).
\end{eqnarray}

\end{document}